\documentclass{siamltex}

\pdfoutput=1
\usepackage[T1]{fontenc}
\usepackage[latin1]{inputenc}
\usepackage{fancybox}

\usepackage{amssymb}
\usepackage{amsmath}
\usepackage{mathpazo}

\usepackage{xcolor}
\usepackage{color}
\usepackage{float}
\usepackage[final]{graphicx}
\usepackage{tikz}
\usetikzlibrary{decorations}
\usetikzlibrary{patterns}

\usepackage{enumerate}

\usepackage[%
  colorlinks%
  ,bookmarks={false}%
  ,pdftitle={paper2}%
  ,pdfsubject={On the Dirichlet Problem for First Order Hyperbolic PDEs on Bounded Domains with Mere Inflow Boundary: Part II Quasi-Linear Equations}%
  ,pdfkeywords={}%
  ,pdfauthor={Thomas März <maerzt@ma.tum.de>}%
  ,plainpages=false
  ]{hyperref}

\input{macros}

\title{On the Dirichlet Problem for First Order Hyperbolic PDEs on Bounded Domains with Mere Inflow Boundary: Part II Quasi-Linear Equations}
\author{Thomas M\"{a}rz\thanks{Zentrum Mathematik,
        Technische Universit\"{a}t M\"{u}nchen,
        Boltzmannstr. 3, 85747 Garching, Germany
        ({\tt maerzt@ma.tum.de}). Manuscript as of \today.}}

\begin{document}

\maketitle

\begin{abstract} 
We study the Dirichlet problem for first order hyperbolic quasi-linear functional PDEs on a simply connected bounded domain of $\R^2$, where the domain has an interior outflow set and a mere inflow boundary. 
While the question of existence of a solution has already been answered in its predecessor,
the present paper discusses the uniqueness and continuous dependence on the coefficients of the PDE.
Under the assumption that the functional dependence is causal, we are able to derive a contraction principle which is the key to proof uniqueness and continuous dependence.
Such a causal functional dependence appears, e.g., in transport based image inpainting. 
\end{abstract}



\section{Introduction}\label{Sect:Intro}
The subject of this paper is hyperbolic quasi-linear functional PDEs of first order in two space dimensions.
We consider the Dirichlet problem on simply connected bounded domains $\Omega$ and bounded functions $u_0$ of bounded variation as boundary data.
The quasi-linear PDE is stated in the space $\BV$ -- the functions of bounded variation -- and formulated in non-conservative form as
\begin{align}\label{eqn:QuasiLinPDE}
	\SP< {c[u](x)\,} , {D u} > &= f[u](x) \cdot \Lm^2  \qquad \mbox{in} \quad \Omega \wo \Sigma \;, \qquad & u|_{\bd \Omega} &= u_0 \;,
\end{align}
where $\SP<\dotarg,\dotarg >$ denotes the scalar product in $\R^2$, $Du$ is the derivative measure of $u \in \BV$ and $\Lm^2$ the Lebesgue measure on $\R^2$.
The dependence on $u$ of the transport field $c$ and the right-hand side $f$ is thereby of a general functional type.
Because the PDE \refEq{eqn:QuasiLinPDE} is hyperbolic, we have to rule out the case of characteristic points.
For this purpose, our requirement is, as in \cite{MyPaper1}, the existence of a time function $T:\Omega \to \R$ which is a global Lyapunov function for the transport field $c$, i.e.,
$c$ is assumed to satisfy the causality condition
\begin{align}\label{eqn:CausalUni}
	\SP< {c[u](x)}, {\nabla T(x)}> \geq \beta \cdot |\nabla T(x)| \;, & \qquad \beta > 0 \;.
\end{align}
In addition, $T$ is zero on the boundary and increases towards its maximal level $\Sigma$, the stop set. Consequently, all boundary points are inflow points and
$\Sigma$ is an interior outflow set.

Our interest is the uniqueness of the solution $u$ and its continuous dependence on the coefficients $c$ and $f$ of PDE \refEq{eqn:QuasiLinPDE}.
In \cite[section 6]{MyPaper1} we have already tackled such quasi-linear equations in following way: 
fix the functional argument of the coefficients by some function $v$ to obtain a linear problem, and the linear theory yields a unique solution $U[v]$ depending on $v$.
Now, every fixed point $u=U[u]$ solves the quasi-linear equation.
We have concluded the existence of a solution as a consequence of the Schauder fixed point theorem, but we have also given an example of non-uniqueness.
The latter example demonstrates that we need further assumptions on the functional argument in order to obtain uniqueness.

Our additional requirement is \emph{functional causality}. With the time function $T$ we have a notion of time on the domain $\Omega$. 
We call the coefficients $c$ and $f$ functionally causal (w.r.t. time $T$) if the values $c[v](x)$ and $f[v](x)$ depend only on those values $v(z)$ which the function $v$ takes on
the subset $\{z \in \Omega : T(z) < T(x)\}$ --  the \emph{past} of $x$.

Functional causality together with Lipschitz continuity of the coefficients (w.r.t. $v$) carry over to the solution operator $U$ and allow us to derive a contraction principle for $U$.
The contraction principle then is the key to show uniqueness and continuous dependence of the fixed point solution $u=U[u]$.

This strategy has been inspired by Kamont's book \cite{Kamont:HFDE}.
In \cite{Kamont:HFDE} the author considers the Cauchy problem for hyperbolic functional PDEs with the functional dependence being causal w.r.t. the physical time
and proofs the local existence of a unique classical solution by employing the Banach fixed point theorem.

Before we get to the point of the theory we should take a moment and think of a possible application. In our case that means transport based \emph{image inpainting}.
The term 
refers to the retouching of undesired or damaged portions of an image.
For this purpose, we have developed, in \cite{FBTM07}, a quasi-linear model for which we can give the following rationale: 
imagine a restorer doing brush strokes in the damaged area $\Omega$. 
Assuming on the one hand  that he only uses color given by the data $u_0$ on $\bd \Omega$
and on the other that brush strokes go along trajectories $x(t)$ -- of a vector field $c$ -- which constantly carry a single color, we end up with the dynamical system
\begin{equation}\label{eqn:IPLinChar}
	\begin{aligned}
		x' &= c(x) \; , & x(0) &= x_0 \in \bd \Omega \;, \\
		u' &= 0 \; , & u(0) &= u_0(x_0) \;,
	\end{aligned}
\end{equation}
which describes exactly the characteristics of a linear advection problem with right-hand side equal to zero. 
Ideally, in order to obtain an aesthetic inpainting, the vector field $c$ would need to reflect the full expertise of our restorer,
which clearly is impossible. But at least the vector field $c$ should be adapted to and hence depend on the image $u$.
This consideration aims for a quasi-linear model of type \refEq{eqn:QuasiLinPDE}.
For our image inpainting model, the dependence of $c[u](x)$ on $u$ is based on an estimation of the tangent vector, which is tangent to the level line of $u$ going through $x$. 
This is because the brush strokes are supposed to continue level lines of $u_0$ which have been interrupted by $\Omega$.
Moreover, the estimation of the tangent only depends on already known information which causes the transport field $c$ to be functionally causal.

In \cite{FBTM07} we discussed the modelling and the construction of a fast algorithm, but we left open the question for the well-posedness of our inpainting model. 
This question is answered positively in \cite{MyDiss}  as an application of the theory presented in this paper.  
\medskip

\paragraph{Outline of the Paper}
In section \ref{Sect:Problem} we review the linear theory and summarize the results of \cite{MyPaper1}. All the requirements of the linear problem are restated here
because they will  be reused in the later sections. The reader already familiar with \cite{MyPaper1} can skip this section.
In section \ref{Sect:QuasiLin} we take up the quasi-linear problem. We collect the requirements on the coefficients of \refEq{eqn:QuasiLinPDE} and
summarize the previous results of \cite{MyPaper1} for the quasi-linear case.
Section \ref{Sect:Unique} is about the uniqueness. After having established the contraction principle mentioned above, we conclude as consequence
the uniqueness in theorem  \ref{Theo:UniquenessQL}.
The subject of section \ref{Sect:ContDep} is the continuous dependence on the coefficients of \refEq{eqn:QuasiLinPDE}.
In the proof of theorem \ref{Theo:ContDepQL} the contraction principle is again the key to the continuous dependence.

\section{Review of the Linear Case}\label{Sect:Problem}
We summarize the theory of the linear problem from \cite{MyPaper1}, because the quasi-linear theory of the later sections is based on it.
We begin with the requirements on the data of the linear problem:
\begin{equation}\label{eqn:LinProblem}
	\begin{aligned}
		 \SP< {c(x)}, {Du} > &= f(x) \cdot \Lm^2 \quad \mbox{ in } \quad \Omega \wo \Sigma  \;,  \\
		 u|_{\bd \Omega} &= u_0   \;, \\ 
		 \SP< {c(x)}, {\nabla T(x)} > & \geq \beta \cdot |\nabla T(x)|  \;.
	\end{aligned}
\end{equation}
The first set of requirements is concerned with the domain $\Omega$, the stop set $\Sigma$ and the time function $T$, while
the second set of requirements is concerned with the transport field $c$, the right-hand side $f$ and the boundary data $u_0$.
\medskip

\begin{req}\label{Req:Domain}(Domains)
	Domains $\Omega \subset \R^2$ are assumed to
	be open, bounded and simply connected
	and to have $\C^1$ boundary.
\end{req}
\medskip

Because of requirement \ref{Req:Domain} the boundary $\bd \Omega$ of a domain is a simple closed $\C^1$ curve. 
By $\gamma:\R \to \bd \Omega$  we denote  a generic regular and periodic parametrization of $\bd \Omega$.
That means $\gamma \in \C^1(\R,\bd \Omega)$ is surjective and $\gamma'(s) \neq 0 \; \forall s \in \R$. 
Furthermore, by $I=\IV[a,b[ \subset \R$ we denote an interval such that $\gamma|_I$ is a generator of $\gamma$.

For our problem time functions are global Lyapunov functions whose range corresponds to a finite time interval.
As time is usually a positive quantity which increases, we define the stop set to be the maximal level of $T$
while in literature the stop set is often the minimal level, see e.g. \cite{Amann:1990:ODE}.
\medskip

\begin{req}\label{Req:TimeFunc}
	(Time functions)
	Time functions are of type \(T: \Omega \to \R \). The upper level-sets of $T$ are denoted by
	\[
		\chi_{T \geq \lambda} := \{x \in \Omega : T(x) \geq \lambda\} \;.
	\]	
	We assume that time functions $T$ satisfy the following conditions:
	\smallskip
	\begin{enumerate}[1.]
		\item $T \in \C(\cl{\Omega})$.
			\smallskip
		\item The boundary of $\Omega$ is the start level: $T|_{\bd \Omega} = 0$. 
			\smallskip
		\item $T$ incorporates a stop set $\Sigma$ with stop time $1$:
			\begin{enumerate}[a)]
				\item $T(x)<1 \Leftrightarrow x \in \cl{\Omega} \wo \Sigma$.
				\item $T|_{\Sigma} = 1$, i.e., $\Sigma$ is the maximal level of $T$. 
			\end{enumerate}
			\smallskip
		\item $T$ increases strictly from $\bd \Omega$ towards $\Sigma$. That means that any upper level-set $\chi_{T \geq \lambda}$ is simply connected and
			\[
				\chi_{T \geq \lambda} = \cl{ \chi_{T > \lambda} } \qquad \forall \; \lambda \in \IV[0,1[ \;.
			\]
		\item Any proper upper level-set is a future domain: for every $\lambda \in [0,1[ $ the set $\chi_{T > \lambda}$ satisfies requirement \ref{Req:Domain}.	
			Furthermore, the field of interior unit normals to the $\lambda$-levels 
			\[
				\chi_{T=\lambda} = \bd \chi_{T > \lambda} \;, \qquad \lambda \in \IV[0,1[ 
			\]
			of $T$ is denoted by $N: \Omega \wo \Sigma \to S^1$. $N$ is required to be continuously differentiable and extendable onto $\bd \Omega$, i.e., $N \in \C^1(\cl{\Omega} \wo \Sigma)$. 
			\smallskip 
		\item[6.*] $T \in \C^2(\cl{\Omega})$, with $\nabla T (x)=0 \Leftrightarrow x \in \Sigma$.
	\end{enumerate}
\end{req}
\medskip

Remark on 6.*: this assumption is in order to ease things in the passages that follow. 
Because of 6.*, we obtain a simple description of the field $N$ on $\Omega \wo \Sigma$ by $	N(x) = \nabla T(x)/|\nabla T(x)|$ and
$N$ is continuously differentiable and extendable onto $\bd \Omega$.
\medskip

In part 3 of requirement \ref{Req:TimeFunc} we have assumed that $T$ features a stop set $\Sigma$. Here we state the geometric properties of allowed stop sets. 
\medskip

\begin{req}\label{Req:StopSet}
	(Stop sets)
	We assume that the stop sets $\Sigma$ satisfy the following conditions:
	\begin{enumerate}[1.]
		\item $\Sigma$ is a closed subset of  $\Omega$.
		\item $\Sigma$ is either an isolated point, or a connected set with tree-like structure (no loops).
		\item If $\Sigma$ is not an isolated point, it consists of finitely many rectifiable $\C^1$ arcs $\Sigma_k$:
			\[
				\Sigma = \bigcup\limits_{k=1}^{n} \Sigma_k \;.
			\]
			The collection $\{ \Sigma_k\}_{ k=1,\ldots,n}$ is assumed to be minimal in the number $n$ of arcs, so $\Sigma$ is decomposed by
			breaking it up at corners and branching points.
			\smallskip
			
			Furthermore, we require for each arc $\Sigma_k$ that its relative interior $\inner{\Sigma}_k$ has a given orientation by a 
			continuous unit normal $n_k : \inner{\Sigma}_k \to S^1$.
	\end{enumerate}
\end{req}
\medskip

In the case in which $\Sigma$ is not just an isolated point, we also need good behavior of the maps $T$ and $N$ close to and on the stop set $\Sigma$. 
For this purpose, we use the following concept of one-sided limits towards $z \in \inner{\Sigma}_k$:

If $x \in B_r(z) \wo \Sigma$ and if $r>0$ is sufficiently small, the projection $p$ of the point $x$ is unique. 
In view of this feature we say a point $x \in B_r(z) \wo \Sigma$
is on the right-hand side or plus side (respectively, on the left-hand or minus side) of $\inner{\Sigma}_k$ if
\begin{equation}
	\frac{x-p}{|x-p|} = +n_k(p)  \qquad \left( \frac{x-p}{|x-p|} = -n_k(p) \right) \;.
\end{equation}
Therewith, a sequence $(x_n)_{n \in \N}$, $x_n \in \Omega \wo \Sigma$, tends to $ z \in \inner{\Sigma}_k $ coming from the plus side 
(respectively, the minus side), in symbols
\begin{equation}
	x_n \to z_+ \qquad \left( x_n \to z_- \right) \;,
\end{equation}
if the sequence converges towards $z$ and almost all elements $x_n$ are on the plus side (respectively, minus side).

With the concept of one-sided limits, the good behavior of the maps $T$ and $N$ is summarized in the following requirement.
\medskip

\begin{req}\label{Req:TimeFuncSigma}	
	(Good behavior at $\Sigma$)
	\begin{enumerate}[1.]
		\item Requirements on $T$:
						
			Let $y \in \Sigma$ and $h \in S^1$ . Let $p=p(y,h)$ be the best possible order for the asymptotic formula
			\[
				T(y + r h) = 1 - \O(r^p) \; ,\quad r \to 0_+ \;.
			\]
			We require that there is a bound $q$ such that $\; \sup\limits_{y \in \Sigma} \sup\limits_{|h|=1} p(y,h) < q$. 
		\item Requirements on $N$:
			\begin{enumerate}[a)]
				\item $N$ has one-sided extensions onto the relatively open components $\inner{\Sigma}_k$ and those extensions are given by $\pm n_k$:
					\begin{align*}
						N^+(y) &\; := \lim\limits_{x \to y_+} N(x) \;, & N^+(y) &= -n_k(y) \;, \\
						N^-(y) &\; :=  \lim\limits_{x \to y_-} N(x) \;, & N^-(y) &= n_k(y) 
					\end{align*}
					for every  $ y \in \inner{\Sigma}_k$. 
					\smallskip
				\item The derivative $DN$ has one-sided extensions onto the relatively open components $\inner{\Sigma}_k$, i.e.,
					\begin{align*}
						(DN)^+(y) &\; := \lim\limits_{x \to y_+} DN(x) \;, &
						(DN)^-(y) &\; :=  \lim\limits_{x \to y_-} DN(x)  
					\end{align*}
					exist for every $y \in \inner{\Sigma}_k$.
					\smallskip
				\item $|DN| \in \L^1(\Omega)$, i.e., poles of $|DN|$ at corner-, branching- and terminal nodes of $\Sigma$ are integrable.
					This feature is assumed to hold in the case in which $\Sigma$ is an isolated point as well.
			\end{enumerate}
	\end{enumerate}
\end{req}
\medskip

So far we have the requirements on domains, stop sets, and time functions. Now we turn to the assumptions on transport fields. 
Here, for a given time function $T$ with stop set $\Sigma$, the causality of the transport field w.r.t. $T$ and its good behavior close to $\Sigma$ are the main concern.
\medskip

\begin{req}\label{Req:TransportField}
	(Transport fields)
	Transport fields of type $c:\Omega \wo \Sigma \to \R^2$ are required to satisfy:
	\medskip
	
	\begin{enumerate}[1.]
		\item $c \in \C^1(\Omega \wo \Sigma)^2$ and $c$ features the following properties:
			\begin{enumerate}[a)]
				\item $c$ and $Dc$ are continuously extendable onto $\bd \Omega $.
				\item If $\Sigma$ is not just an isolated point, then $c$ and $Dc$ have one-sided limits on the relatively open $C^1$ arcs
					$\inner{\Sigma}_k$ of $\Sigma$:
					\begin{align*}
						c^+(y) = & \; \lim\limits_{x \to y_+} c(x) \quad \mbox{and} \quad c^-(y) = \lim\limits_{x \to y_-} c(x) \; , \\
						(Dc)^+(y) = & \; \lim\limits_{x \to y_+} Dc(x) \quad \mbox{and} \quad (Dc)^-(y) = \lim\limits_{x \to y_-} Dc(x) \; ,
					\end{align*}
					for every $y \in \inner{\Sigma}_k$.
			\end{enumerate}
			\medskip
		\item Unit speed and causality condition:
			\begin{enumerate}[a)]
				\item $|c| = 1$.
				\item There is a lower bound $\beta > 0$ such that
					\[ 
						\beta \leq \SP< {c(x)}, {N(x)} > \leq 1 \qquad \forall x \in \cl{\Omega} \wo \Sigma \:.
					\]
				\item Conditions a) and b) hold for the one-sided limits as well, i.e.,
					\[
						|c^+(y)| = |c^-(y)| = 1
					\]
					and
					\[ 
						\beta \leq \SP< {c^+(y)}, {N^+(y)} > \leq 1 \quad , \quad \beta \leq \SP< {c^-(y)}, {N^-(y)} > \leq 1 \; , 
					\] 
					whenever $y$ belongs to some  $\inner{\Sigma}_k$.
			\end{enumerate}
			\medskip
		\item Let $z_k$, $k \in \{1, \ldots, m\}$ denote the terminal-, branching- and kink nodes of $\Sigma$. 			
			For every $\varepsilon > 0$ such that each disk $B_{\varepsilon}(z_k)$
			is compactly contained in $\Omega$, we define the set
			\[
				V_{\varepsilon} := \Sigma \; \cup \;  \bigcup\limits_{k=1}^m \cl{ B_{\varepsilon}(z_k) } \;.
			\]
			\begin{enumerate}[a)]
				\item For every admissible $\varepsilon > 0$, there is a bound $M_{\varepsilon}$ such that
					\[
						|Dc(x)| \leq M_{\varepsilon} \quad , \quad \forall \; x \in \Omega \wo V_{\varepsilon}
					\]
				\item $|Dc| \in \L^1(\Omega)$, i.e.,  poles of $|Dc|$ at $z_k$, $k \in \{1, \ldots, n\}$ are integrable.		
			\end{enumerate}
			\medskip
	\end{enumerate}
\end{req}
\medskip

Finally we need a right-hand side and boundary data.
\medskip

\begin{req}\label{Req:RHSandData}
	(Right-hand sides and boundary data)
	We assume right-hand sides $f: \Omega \to \R$ and Dirichlet boundary data $u_0 : \bd \Omega \to \R$ to be $f \in C^1(\cl{\Omega})$ and $u_0 \in \BV(\bd \Omega)$.
\end{req}
\smallskip

Remark: with $u_0 \in \BV(\bd \Omega)$ we mean that for every regular periodic parametri\-zation $\gamma$ of $\bd \Omega$, the pull-back $\gamma^* u_0 = u_0 \circ \gamma$ is a periodic $\BV$ function.
Moreover, because $\bd \Omega$ is one-dimensional, $\gamma^* u_0$ is $\BV$ function of one variable. Therefor the boundary data is essentially bounded.

In the following sections we assume, if not explicitly stated otherwise, that domains $\Omega$, stop sets $\Sigma$, time functions $T$, and boundary data $u_0$
satisfy the requirements above. For the quasi-linear case we will add up the list of requirements on transport fields $c$ and right-hand sides $f$ later on.
\medskip

\paragraph{Results of \cite{MyPaper1}}
Problem \refEq{eqn:LinProblem} has a unique solution $u$ in $\BV(\Omega)$ which depends continuously on all the data of the problem.
Employing the method of characteristics the solution is given by
\begin{equation}\label{eqn:SolOrig} 
	u(x) = u_0( \eta(T_0(x),x)) + \int\limits_0^{T_0(x)} f_0 \circ \eta(\tau ,x) \: d\tau \;. 
\end{equation}
Here, $\eta$ denotes the backward characteristics which are the solution of
\begin{align}\label{eqn:BCharIVP}
	\eta'(\dotarg,x) &= -c_0 \circ \eta(\dotarg,x) \;, & \eta(0,x) &= x \in \Omega \wo \Sigma \;.
\end{align} 
The forward characteristics -- denoted by $\xi$ -- solve
\begin{align}\label{eqn:FCharIVP}
	\xi'(\dotarg,s) &= c_0 \circ \xi(\dotarg,s) \;, & \xi(0,s) &= \gamma(s) \in \bd \Omega \;.
\end{align} 
According to \cite[section 2]{MyPaper1},
$T_0=1-(1-T(x))^{1/q}$ denotes the transformed version of $T$ whose gradient $\nabla T_0$ blows up at $\Sigma$. $T_0$ is equivalent to $T$ for it has the same level lines as $T$ (level lines get only new names).
With $T_0$ we also transform the transport field and the right-hand side by $c_0 = c / \SP< {c}, {\nabla T_0}>$ and $f_0 = f / \SP< {c}, {\nabla T_0}>$. 
The forward characteristics w.r.t. $c_0$ have then the useful property $T_0( \xi(t,s) ) = t $, i.e.,
$T_0$ denotes the time of the characteristics. Moreover, $T_0$ features the properties $T_0 |_{\Sigma} =1$, $T_0 |_{\bd \Omega} =0$, and $|\nabla T_0(x)| \geq m_0 > 0$ which
imply the bound $|c_0| \leq 1/m_0 \beta$, a time range of $t \in [0,1]$ for $\xi(\dotarg,s)$ , and a bound of $1/m_0 \beta$ on the arc length of $\xi(\dotarg,s)$.

By $v(t,s) = u \circ \xi(t,s)$ the solution $u$ is represented in characteristic variables $(t,s)$ as
\begin{equation}\label{eqn:SolChar}
	v(t,s) = \gamma^*u_0(s) + \int\limits_0^t f_0 \circ \xi(\tau,s) \: d\tau \;. 
\end{equation}

Theorem 4.1 of \cite{MyPaper1} shows that $u \in \BV(\Omega)$ and not only $u \in \BV(\Omega \wo \Sigma)$. 
The continuous dependences of the solution $u$ on $c$ and $T$ are given by theorem 5.1 of \cite{MyPaper1},
while the continuous dependences on $u_0$ and $f$ are part of theorem 4.1 .

This completes the summary of the linear theory and we now turn to the quasi-linear problem.

\section{The Quasi-Linear Problem}\label{Sect:QuasiLin}
In \cite[section 6]{MyPaper1} we have already shown that the quasi-linear problem 
\begin{equation}\label{eqn:QuasiLinProblem}
	\begin{aligned}
		 \SP< {c[u](x)}, {Du} > &= f[u](x) \cdot \Lm^2 \quad \mbox{ in } \quad \Omega \wo \Sigma  \;, \\
		 u|_{\bd \Omega} &= u_0   \;, \\  
		 \SP< {c[u](x)}, {\nabla T(x)} > & \geq \beta \cdot |\nabla T(x)|   \;,
	\end{aligned}
\end{equation}
where the dependencies of $c$ and $f$ on $u$ are of a general functional type, admits a solution in $\BV(\Omega)$.
We achieved this result as a consequence of the Schauder fixed point theorem applied to the operator $U[v]$, which denotes the solution operator to the family of linear problems
\begin{equation}\label{eqn:LinOfQuasiLin}
	\begin{aligned}
		 \SP< {c[v](x)}, {Du} > &= f[v](x) \cdot \Lm^2 \quad \mbox{ in } \quad \Omega \wo \Sigma  \;, \quad v \in \FX \subset \BV(\Omega)\\
		 u|_{\bd \Omega} &= u_0   \;, \\  
		 \SP< {c[v](x)}, {\nabla T(x)} > & \geq \beta \cdot |\nabla T(x)|   \;.
	\end{aligned}
\end{equation}
In addition, we have given an example of non-uniqueness which tells us that we need stronger requirements on $c$ and $f$
in order to obtain uniqueness.

The requirements used here which will suffice to ensure uniqueness and beyond that continuous dependence are summarized in the following.
A special ingredient is causality w.r.t. the functional argument $v$, that is the dependence on $v$ is of Volterra type w.r.t. the notion of time induced by our time function $T$.
\medskip

\begin{defi}
	(Functional causality) 
	Let $T(x)$ be the time of the point $x \in \Omega$.  
	The set $\Omega_{T(x)} := \{z \in \Omega: T(z) < T(x)\}$,
	denotes the past of the point $x$ (w.r.t. the time function $T$).
	Let $\FD(\Omega)$ and $\FR(\Omega)$ be function spaces defined on $\Omega$ and let $f$ be an operator of type
	\begin{align*}
		f:\FD(\Omega) & \to \FR(\Omega) \;, \quad \mbox{with} \; & f[\dotarg](x):\FD(\Omega) & \to \R^d \;, \quad  x \in \Omega \;.
	\end{align*}	
	We say that $f$ is functionally causal (w.r.t. time $T$) if 
	\[
		f[v](x) = f[ \; v \cdot \1_{\Omega_{T(x)}}\;] (x) \;.
	\]
\end{defi}
\smallskip

This definition says that the value $f[v](x)$ depends only those values $v(z)$ which $v$ attains on the past $\Omega_{T(x)}$.
Beyond that the functional causality implies the following domain restriction feature. Let $\lambda$ be in the range of $T$.
Then, for $x \in \Omega_{\lambda} = \{z \in \Omega: T(z) < \lambda\}$ and $v \in \FD(\Omega)$, the inclusion $\Omega_{T(x)} \subset \Omega_{\lambda}$ implies
\[
	f[v](x) = f[ \; v \cdot \1_{\Omega_{T(x)}} \;] (x) = f[ \; v \cdot \1_{\Omega_{\lambda}} \;] (x) \;.
\]
Hence, the domain restriction (onto $\Omega_{\lambda}$ ) by $f_\lambda:\FD(\Omega_{\lambda})  \to \FR(\Omega_{\lambda})$, $f_\lambda [v] = f[v]|_{\Omega_\lambda}$ is well-defined.
\medskip

\begin{req}\label{Req:TransportFieldFunc}
	(Transport fields)
	Regarding the functional argument, transport fields are maps of type
	\begin{align*}
		c: \L^1(\Omega) & \to \C^1(\Omega \wo \Sigma)^2 \;, & \mbox{with} & & c[\dotarg](x): \L^1(\Omega) &\to \R^2 \;,
	\end{align*}
	and we assume them to satisfy:
	\smallskip
	
	\begin{enumerate}[1.]
		\item $c$ is functionally causal.
		\item For fixed $v \in \L^1(\Omega)$ the function $c[v]: \Omega \wo \Sigma \to \R^2 $ is a transport field according to requirement \ref{Req:TransportField}.
			\smallskip
		\item Uniformity of the unit speed and causality (w.r.t. $x$) condition:
			\smallskip
			\begin{enumerate}[a)]
				\item $|c[v](x)| = 1$ for all $x \in \Omega \wo \Sigma$ and for all $v \in \L^1(\Omega)$.
					\smallskip
				\item There is a uniform lower bound $\beta > 0$ such that
					\begin{align*}
						\beta \leq \SP< {c[v](x)}, {N(x)} >  &\leq 1 \qquad \forall \; x \in \cl{\Omega} \wo \Sigma 
						& \mbox{and}\quad &\; \forall \;  v \in \L^1(\Omega) \;.
					\end{align*}
				\item Both conditions hold for the one-sided limits of $c[v]$ on the relatively open $\C^1$ arcs $\inner{\Sigma}_k$ of $\Sigma$.
			\end{enumerate}
			\smallskip
		\item Bounds and continuity:
			\smallskip
			\begin{enumerate}[a)]
				\item The map $D_x c: \L^1(\Omega) \to \C(\Omega \wo \Sigma)^{2 \times 2}$ -- the derivative of $c[v]$ w.r.t. the variable $x$ -- 
					is $\L^1$ bounded by 
					\[
						\|D_x c[v]\|_{L^1(\Omega)} < M_1 \qquad \forall \; v \in \L^1(\Omega) \;.
					\] 
				\item the bound on $|D_x c[v](x)|$, from requirement \ref{Req:TransportField} 3a),
					\[
						|D_x c[v](x)| \leq M_{\varepsilon} \quad , \quad \forall \; x \in \Omega \wo V_{\varepsilon}
					\]
					holds uniformly for all $v \in \L^1(\Omega)$.
				\item the map $c$ is Lipschitz in the following manner:
					\[
						\|c[v] - c[w]\|_{\infty} \leq L_1 \cdot \|v-w\|_{L^1(\Omega)} \;.
					\]
			\end{enumerate}
	\end{enumerate}
\end{req}
\medskip

\begin{req}\label{Req:RHSFunc}
	(Right-hand sides)
	Regarding the functional argument, right-hand sides are maps of type
	\begin{align*}
		f: \L^1(\Omega) & \to \C^1(\cl{\Omega}) \;, & \mbox{with} & & f[\dotarg](x): \L^1(\Omega) &\to \R \;,
	\end{align*}
	and we assume them to satisfy:
	\smallskip
	\begin{enumerate}[1.]
	\item $f$ is functionally causal.
	\item Bounds and continuity:
		\begin{enumerate}[a)]
			\item The map $f$ is bounded by 
				\[
					 \|f[v]\|_{\infty} \leq M_2 \qquad \forall \; v \in \L^1(\Omega) \;.
				\] 
			\item The map $\nabla_x f: \L^1(\Omega) \to \C(\cl{\Omega})^2$ -- the derivative of $f[v]$ w.r.t. the variable $x$ -- is bounded by
				\[
					 \|\nabla_x f[v]\|_{\infty} \leq M_3 \qquad \forall \; v \in \L^1(\Omega) \;.
				\]
			\item the map $f$ is Lipschitz in the following manner:
				\[
					\|f[v] - f[w]\|_{\infty} \leq L_2 \cdot \|v-w\|_{L^1(\Omega)} \;.
				\]
		\end{enumerate}
	\end{enumerate}
\end{req}
\medskip

Finally, we define, as in \cite{MyPaper1}, the subsets of $\BV(\bd \Omega)$ and $\BV(\Omega)$ with which we will work later on.
\medskip

\begin{defi}\label{Def:FuncDomains}
	Let $M_1$, $M_2$, $M_3$ be the bounds from the requirements stated above. 
	\begin{enumerate}[a)]
		\item We denote by
			\[
				\FB = \FB(\bd \Omega) := \{v \in \BV(\bd \Omega) : \|v\|_{\L^{\infty}(\bd \Omega)} \leq M_4 \;,\; |Dv| \leq M_5\}
			\]
			the set of boundary functions.
		\item Let the constants $M_{*} \in \R$ and $M_{**} \in \R$ be given by
			\begin{equation}\label{eqn:SelfMapBound}
			\begin{aligned}
				M_{*}  := &\left(M_4 + \frac{M_2}{\beta \cdot m_0}\right) \cdot \Lm^2(\Omega) \;, \\
				M_{**} := & \; 2 \cdot \left(M_4 + \frac{M_2}{\beta \cdot m_0}\right) \cdot \Hm^1(\Sigma) + \frac{M_5}{\beta \cdot m_0}  + \left( \frac{M_2}{\beta} + \frac{M_3}{\beta^2 \cdot m_0}\right) \cdot \Lm^2(\Omega) \\
				& + \frac{M_2}{\beta^3 \cdot m_0} \cdot \left( M_1 + \|DN\|_{L^1(\Omega)}\right) \;.
			\end{aligned}
			\end{equation}
			We set
			\[
				\FX = \FX(\Omega) := \{ v \in \BV(\Omega) : \|v\|_{L^1(\Omega)} \leq M_{*} \; ,\; |Dv|(\Omega) \leq M_{**}\} \;.
			\]
	\end{enumerate}
\end{defi}
\medskip

When speaking about transport fields $c$ and right-hand sides $f$ in the following, we mean transport fields and right-hand sides according to requirements \ref{Req:TransportFieldFunc} and \ref{Req:RHSFunc}.

In \cite{MyPaper1} we have shown that the solution operator $U:\FX \to \FX$ is a self-mapping. Here, $U$ inherits additionally the functional causality of $c$ and $f$.
\medskip

\begin{lem}\label{Lem:VolterraU}
	The solution operator $U:\FX \to \FX$ of problem \refEq{eqn:LinOfQuasiLin} is functionally causal.
\end{lem}
\smallskip

\begin{proof}
	Let $v \in \FX$ be arbitrary but fixed. Then, according to equation \refEq{eqn:SolOrig}, $U[v](x)$ is given by
	\[
		U[v](x) = u_0( \eta[v](T_0(x),x)) + \int\limits_0^{T_0(x)} f_0[v] \circ \eta[v](\tau ,x) \: d\tau \;,
	\]
	where the backward characteristic $\eta[v](\dotarg,x)$ is the solution of
	\begin{align*}
		y' &= -c_0[v](y) \;, & y(0) &= x \;,
	\end{align*}
	while $c_0$ and $f_0$ are given by
	\begin{align*}
		c_0[v](x) &= \frac{c[v](x)}{\SP< {c[v](x)}, {\nabla T_0(x)}>} \;,& f_0[v](x) &= \frac{f[v](x)}{\SP< {c[v](x)}, {\nabla T_0(x)}>} \;.
	\end{align*}
	Clearly, $c_0$ and $f_0$ are the same way functionally causal as $c$ and $f$.
	
	Since $T_0$ has the same level-sets as $T$ while merely the names of level lines have changed, we refer to $T_0$ in order to denote the past
	$\Omega_{T_0(x)} := \{z \in \Omega : T_0(z) < T_0(x)\}$.
	
	For every $t \in \IV]0,{T_0(x)}[$ we know that $\eta[v](t,x) \in \Omega_{T_0(x)}$. The latter is a consequence of the uniform causality w.r.t. $x$ (requirement \ref{Req:TransportFieldFunc} part 3b)).
	Hence,  $\eta[v](\dotarg,x)$ depends merely on the restriction of the transport field $c_0[v]$ onto the past $\Omega_{T_0(x)}$. 
	By the functional causality, $c_0[v]|_{ \Omega_{T_0(x)} }$ w.r.t. $v$ depends only on $v \cdot \1_{\Omega_{T_0(x)}}$. 
	Consequently, for $t \in \IV]0,{T_0(x)}[$ the expression $\eta[v](t,x)$ w.r.t. $v$ also depends only on $v \cdot \1_{\Omega_{T_0(x)}}$.
	Thus, by the representation of $U[v](x)$ above, it is obvious that $U$ is functionally causal:
	\[
		U[v](x) = U[\; v \cdot \1_{\Omega_{T_0(x)}} \;](x) \;. 
	\]
	\hfill
\end{proof}
\medskip

Using the domain restriction feature, as discussed above, we define $\FX_{\lambda}$: 
\medskip

\begin{defi}
	Let $\lambda$ be in the range of $T$. We denote by
	$\FX_{\lambda} = \FX(\Omega_{\lambda}) := \{ v|_{\Omega_{\lambda}} : v \in \FX \}$
	the domain-restricted version of $\FX$ from definition \ref{Def:FuncDomains}.
\end{defi}
\medskip

In the following, we consider the domain restricted quasi-linear problem on $\Omega_\lambda$:
\begin{equation}\label{eqn:QuasiLinProblemRest}
	\begin{aligned}
		 \SP< {c[u](x)}, {Du} > &= f[u](x) \cdot \Lm^2 \quad \mbox{ in } \quad \Omega_\lambda  \;, \\
		 u|_{\bd \Omega} &= u_0   \;, \\  
		 \SP< {c[u](x)}, {\nabla T(x)} > & \geq \beta \cdot |\nabla T(x)|   \;,
	\end{aligned}
\end{equation}
The question for existence has already been answered in \cite[section 6]{MyPaper1}: 
by lemma \ref{Lem:VolterraU} the operator $U:\FX \to \FX$ is functionally causal. Hence, the domain restricted version $U_\lambda:\FX_{\lambda} \to \FX_{\lambda}$
is well-defined and $U_\lambda[v]$ solves the domain restricted linear problem
\begin{equation}\label{eqn:LinOfQuasiLinRest}
	\begin{aligned}
		 \SP< {c[v](x)}, {Du} > &= f[v](x) \cdot \Lm^2 \quad \mbox{ in } \quad \Omega_\lambda  \;, \quad v \in \FX_\lambda \\
		 u|_{\bd \Omega} &= u_0   \;, \\  
		 \SP< {c[v](x)}, {\nabla T(x)} > & \geq \beta \cdot |\nabla T(x)|   \;.
	\end{aligned}
\end{equation}
Since $U_\lambda:\FX_{\lambda} \to \FX_{\lambda}$ is continuous (w.r.t. the $\BV$ weak* topology) and $\FX_\lambda$ is non-empty, convex, and compact 
(see \cite[section 6]{MyPaper1} for both statements), we get a solution of \refEq{eqn:QuasiLinProblemRest} as a fixed point of $U_{\lambda}$
by the Schauder-Tychonoff theorem.
And vice versa, every fixed point $u$ of the original operator $U:\FX \to \FX$ after restriction $u|_{\Omega_{\lambda}}$ belongs to $\FX_{\lambda}$ and is a fixed point of $U_\lambda:\FX_{\lambda} \to \FX_{\lambda}$.

In the following, we show that for every $\lambda$ in the range of $T$ the fixed point of $U_\lambda$ is unique.

\section{Uniqueness of the Fixed Point}\label{Sect:Unique}
The subject of this section is to show that the solution of problem \refEq{eqn:QuasiLinProblem} is unique or rather that $U$ has a unique fixed point.
For this purpose we will show that, for any choice of $0 < \lambda < 1$, the domain restricted operator $U_\lambda:\FX_{\lambda} \to \FX_{\lambda}$ is Lipschitz.
Moreover, we will see that $U_\lambda$ is in fact contractive for a suitable choice of $\lambda$. The contractiveness will then be the key
to the uniqueness of the fixed point $u=U[u]$.

In order to estimate the difference $U[v_1] - U[v_2]$ we prepare by setting up a PDE which is satisfied by the difference.
For the purpose of abbreviation, we set
\begin{align*}
	c_1 &:= c[v_1]  \;, & c_2 &:= c[v_2] \;, \quad \mbox{and} \quad &
	f_1 &:= f[v_1]  \;, & f_2 &:= f[v_2] \;.
\end{align*}
Let $u_1$ and $u_2$ respectively denote the solutions of the two linear problems
\begin{align*}
	\SP< {c_1(x)}, {Du} > &= f_1(x) \cdot \Lm^2 \quad \mbox{ in } \quad \Omega_{\lambda} \;,  && u|_{\bd \Omega} = u_{0,1} \; ,\\
	\mbox{and} \quad \SP< {c_2(x)}, {Du} > &= f_2(x) \cdot \Lm^2  \quad \mbox{ in } \quad \Omega_{\lambda} \;,  && u|_{\bd \Omega} = u_{0,2} \; .
\end{align*}
As in the proof of lemma \ref{Lem:VolterraU} we refer to the transformed time $T_0$ instead of $T$ and denote by $\Omega_{\lambda}$ the 
lower level-set of $T_0$
\begin{align*}
	\Omega_{\lambda}  &:= \{z \in \Omega : T_0(z) < \lambda\} \;, & \lambda & \; \in \IV]0,1[ \;.
\end{align*}
For the first considerations we use different boundary data. When setting $u_{0,1} = u_{0,2} = u_0$ later on, we will obtain the relations
\begin{equation}\label{eqn:PrepareU}
	\begin{aligned}
		u_1 &= U[v_1]  \;, & u_2 &= U[v_2] \;.
	\end{aligned}
\end{equation}
Let $w$ denote the difference  $w := u_1-u_2$.  After having subtracted the problems from each other, the difference $w$ must satisfy the linear problem
\begin{align*}
	\SP< {c_1(x)}, {Dw} > &= (f_1(x)-f_2(x)) \cdot \Lm^2  - \SP< {c_1(x)-c_2(x)},{Du_2} > \quad \mbox{ in } \quad \Omega_{\lambda} \;, \\
	w|_{\bd \Omega}  &= w_0 \;,
\end{align*}
with boundary data $w_0 = u_{0,1}-u_{0,2}$.

By the same argumentation as used in \cite[section 5]{MyPaper1} we see that $w$ is the unique solution of this PDE. But, in order to solve for $w$, we cannot directly
apply the method of characteristics as in \cite[section 4]{MyPaper1}, since the right-hand side is not an absolutely continuous measure.

Instead we approximate the right-hand side by absolutely continuous measures. 
Since $u_2 \in \BV(\Omega)$, there exists a sequence $(u_{2,n})_{n \in \N} $ of $\C^{\infty}(\Omega)$ functions
which converges strictly to $u_{2,n}$, i.e.,
\[
	\|u_2-u_{2,n}\|_{\L^1(\Omega)} \to 0 \quad \mbox{and} \quad ||Du_2|(\Omega) - |Du_{2,n}|(\Omega)| \to 0 \;.
\]
For this statement see e.g. \cite[theorem 3.9]{AmbrosioBV}. 
Moreover, we have $Du_{2,n} = \nabla u_{2,n}(x) \cdot \Lm^2$ and $|Du_{2,n}|(\Omega) = \|\nabla u_{2,n}\|_{L^1(\Omega)}$. 
\smallskip

Using such a sequence we obtain an approximate problem
\begin{equation}\label{eqn:PDEDiffApprox}
	\begin{aligned}
		\SP< {c_1(x)}, {Dw} > &= \left( f_1(x)-f_2(x)   - \SP< {c_1(x)-c_2(x)},{\nabla u_{2,n}(x)} > \right) \cdot \Lm^2 \quad \mbox{ in } \quad \Omega_{\lambda} \;, \\
		w|_{\bd \Omega}  &= w_0 \;,
	\end{aligned}
\end{equation}
with a sequence of solutions $w_n$ which we can construct using the method of characteristics. 
By scaling PDE \refEq{eqn:PDEDiffApprox} by the factor $1/\left<c_1,\nabla T_0\right>$, the family of  forward characteristics $\xi(\dotarg,s)$ is then given by the IVP
\begin{align*}
	y' &= c_{1,0} (y) \;,& y(0) &= \gamma(s) \;.
\end{align*}
Here, we set
\begin{align*}
	c_{1,0} :=& \; \frac{c_1}{\left<c_1,\nabla T_0\right>} \;, & f_{1,0} :=&\; \frac{f_1}{\left<c_1,\nabla T_0\right>}  \;,   \quad \mbox{and} \; &
	c_2^0 :=& \; \frac{c_2}{\left<c_1,\nabla T_0\right>} \;, & f_2^0 :=& \; \frac{f_2}{\left<c_1,\nabla T_0\right>}  \;. 
\end{align*}
(Note: if we were to be consistent, we would set $f_{2,0} := f_2/\left<c_2,\nabla T_0\right>$, which differs from $f_2^0$.)
According to equation \refEq{eqn:SolChar}, we obtain $w_n$ in characteristic variables by
\begin{equation}\label{eqn:RepWn}
	w_n \circ \xi(t,s) = \gamma^*w_0(s) + \int\limits_0^t \left(f_{1,0}-f_2^0  - \SP< {c_{1,0}-c_2^0},{\nabla u_{2,n}} >\right) \circ \xi(\tau,s) \; d\tau \;.
\end{equation}
The consideration of the sequence $w_n$ will not be of any use if $w_n$ does not tend to $w$ in an appropriate fashion.
We will show the desired convergence in lemma \ref{Lem:ConvWWn}. But first, we rewrite $w_n \circ \xi$ in order to get a more convenient representation of the difference $(w_n - w) \circ \xi(t,s)$.

Because the PDE for $u_1$ has the same transport field $c_1$, we write $u_1$ in characteristic variables as
\begin{align*}
	u_1 \circ \xi(t,s) &= \gamma^*u_{0,1} + \int\limits_0^t f_{1,0} \circ \xi(\tau,s) \; d\tau \;.
\end{align*}
And because $u_{2,n}$ is smooth, we have
\begin{align*}
	u_{2,n} \circ \xi(t,s) - u_{2,n} \circ \gamma(s) &=  \int\limits_0^t \SP< {c_{1,0}},{\nabla u_{2,n}} > \circ \xi(\tau,s) \; d\tau \;.
\end{align*}
Together with $w_0 = u_{0,1}-u_{0,2}$, the last two observations imply
\begin{align*}
	w_n \circ \xi(t,s) = & \; u_1 \circ \xi(t,s)- u_{2,n} \circ \xi(t,s) + u_{2,n} \circ \gamma(s) - \gamma^*u_{0,2}(s) \\
	& \; + \int\limits_0^t \left( \SP< {c_2^0},{\nabla u_{2,n}} > - f_2^0 \right) \circ \xi(\tau,s) \; d\tau \;.
\end{align*}
Subtracting $w = u_1 - u_2$ finally, we end up with
\begin{equation}\label{eqn:ReqDWWn}
	\begin{aligned}
		(w_n - w) \circ \xi(t,s) = & \; (u_2 - u_{2,n}) \circ \xi(t,s) \; + (\gamma^*u_{2,n}(s) - \gamma^*u_{0,2}(s)) \\
		& \; + \int\limits_0^t \left( \SP< {c_2^0},{\nabla u_{2,n}} > - f_2^0\right) \circ \xi(\tau,s) \; d\tau \;.
	\end{aligned} 
\end{equation}

As a second step of preparation, we show that requirement \ref{Req:TransportFieldFunc} part 4 b) implies uniform bounds on the determinant of $D\xi$.
\medskip

\begin{lem}\label{Lem:UniformDetBound}
	Let $c$ be a transport field. 
	For fixed $v \in \L^1(\Omega)$ let $\xi[v]$  
	denote the general solution of the IVP
	\begin{align}\label{eqn:FuncIVP}
		y' &= c_0[v](y) \;, & y(0) &= \gamma(s) \;, & \mbox{with} \quad c_0[v](y)&=\frac{c[v](y)}{\SP<{c[v](y)}, {\nabla T_0(y)} >} \;.
	\end{align}
	Then,
	\begin{enumerate}[a)]
	\item $\xi[v]: \IV]0,1[ \times \IV]a,b[ \to \Omega \wo (S \cup \Sigma)$ is a diffeomorphism, where
	$S := \{\xi[v](t,a): t \in ]0,1[\}$. Here, the domain of $s \in \IV]a,b[$ refers to the domain of the boundary parametrization $\gamma$, while
	the time domain $t \in \IV]0,1[$ refers to the range of $T_0$.
	\medskip
	\item for $0 < \lambda < 1$ there are bounds $k_{\lambda}$ and $K_{\lambda}$ such that
	\begin{align*}
		0 < k_{\lambda} & \leq \det D\xi[v](t,s) \leq K_{\lambda} \; & \forall \; (t,s) &  \in \IV]0,\lambda[ \times \IV]a,b[ \;.
	\end{align*}
	The bounds $k_{\lambda}$ and $K_{\lambda}$ depend only on $\lambda$, but not on $v$.
	Moreover, $k_{\lambda}$ decreases, while $K_{\lambda}$ increases monotonically with $\lambda$.
	\end{enumerate}
\end{lem}
\medskip

\begin{proof} Statement a) is taken from \cite[corollary 3.3]{MyPaper1}, so only statement b) needs a proof.
	On the restricted domain $\Omega_\lambda$ there is -- by requirement \ref{Req:TransportFieldFunc} part 4 b) -- a uniform bound on the derivative of $c[v]$
	\[
		|D_x c[v](x)| \leq M_{\lambda} \quad , \quad \forall \; x \in \Omega_\lambda \quad , \quad \forall \; v \in \L^1(\Omega) \;.
	\]
	A similar bound -- which we also call $M_{\lambda}$ -- will hold for the derivative $D_x c_0[v](x)$ of the transformed transport field $c_0[v]$.
	
	Because $\xi[v]$ solves the IVP \refEq{eqn:FuncIVP} and $|c_0[v](y)| \leq 1/m_0 \beta$, we have a first estimate
	\[
		\det D\xi[v] \leq |\pd_t \xi[v]| \cdot |\pd_s \xi[v]| \leq \frac{|\pd_s \xi[v]|}{\beta \cdot m_0} \;.
	\]
	The derivative $\pd_s \xi[v]$ is the solution of the variational equation below
	\begin{align*}
		\pd_s \xi[v]' &= D_x c_0[v] \circ \xi[v] \cdot \pd_s \xi[v]  \; & y(0) &= \gamma'(s) \;.
	\end{align*}
	For every $t \in \IV]0,\lambda[$ the point $\xi(t,s)$ belongs to $\Omega_\lambda$. Thus, we estimate
	\[
		|\pd_s \xi[v]|(t,s) \leq \|\gamma'\|_{\infty} + \int\limits_{0}^t M_{\lambda} \cdot |\pd_s \xi[v]|(\tau,s) \; d\tau \;,
	\]
	and an application of Gronwall's lemma (see e.g. \cite{Walter:DiffInEq}) yields
	\[
		|\pd_s \xi[v]|(t,s) \leq \|\gamma'\|_{\infty}  \exp\left(\lambda \cdot M_{\lambda}\right) \;.
	\]
	Hence, we obtain the upper bound $K_{\lambda}$
	\[
		\det D\xi[v] (t,s) \leq \frac{\|\gamma'\|_{\infty}  \exp\left(\lambda \cdot M_{\lambda}\right) }{\beta \cdot m_0} =: K_{\lambda}
	\]
	on $\IV]0,\lambda[ \times \IV]a,b[ $.
	\smallskip
	
	For the lower bound we consider the inverse map $\xi[v]^{-1}(x)$ for $x \in \Omega_{\lambda} \wo S$. According to \cite[corollary 3.3]{MyPaper1}, the inverse map is given by
	\[
		\xi[v]^{-1}(x)= (T_0(x) , s[v](x)) \;, \quad \mbox{with} \quad s[v](x) = \gamma^{-1}( \eta[v](T_0(x),x)) \;.
	\]
	Therein, $\eta[v](\dotarg,x)$ denotes the backward characteristics given as the solution of
	\begin{align*}
		y' &= -c_0[v](y) \; & y(0) &= x  \in \Omega_{\lambda} \wo S \;.
	\end{align*}
	Now, the determinant of $D_x \xi[v]^{-1}$ is bounded by
	\[
		\det D_x \xi[v]^{-1}(x) \leq |\nabla T_0(x)| \cdot |\nabla_x s[v](x)|
	\]
	and $\nabla_x s[v](x)^T$ is given by the expression
	\[
		 (\gamma^{-1})'(\eta)^T \cdot \left( \pd_t \eta[v](T_0(x),x) \cdot \nabla T_0(x)^T + D_x \eta[v](t,x)|_{t=T_0(x)} \right) \;.
	\]
	The derivative $|D_x \eta[v](t,x)|$ can be estimated in the same way as $|\pd_s \xi(t,s)|$. And because $x \in \Omega_{\lambda}$, we obtain
	\[
		|D_x \eta[v](t,x)| \leq \exp\left(\lambda \cdot M_{\lambda}\right) \;.
	\]
	Using the latter result to establish a bound on $|\nabla_x s[v](x)|$, we end up with
	\[
		\det D_x \xi[v]^{-1}(x) \leq \frac{\|\nabla T_0\|_{\L^{\infty}(\Omega_{\lambda})}}{\min_{s \in [a,b]} |\gamma'(s)|} \cdot 
		\left( \frac{\|\nabla T_0\|_{\L^{\infty}(\Omega_{\lambda})}}{\beta \cdot m_0} + \exp\left(\lambda \cdot M_{\lambda}\right) \right) \;,
	\]
	for $x \in \Omega_{\lambda} \wo S$.
	Finally, we define $1/k_{\lambda}$ equal to the right-hand side of the last inequality and get the desired lower bound  
	\[
		\frac{1}{\det D\xi[v](t,s)} = \det \left( (D\xi[v](t,s))^{-1} \right)= \det D_x \xi[v]^{-1}(x)|_{x = \xi(t,s)} \leq \frac{1}{k_{\lambda}} \;,
	\]
	for $(t,s) \in \IV]0,\lambda[ \times \IV]a,b[ $, since in this case $\xi(t,s) \in \Omega_{\lambda} \wo S$.
		
	Both bounds do not depend on the choice of $v \in \L^1(\Omega)$ and 
	the monotonicity properties of $k_{\lambda}$ and $K_{\lambda}$ as functions of $\lambda$ are obvious. \hfill
\end{proof}
\medskip

In the next lemma we turn to the approximation of $w$ by $w_n$.
\medskip

\begin{lem}\label{Lem:ConvWWn}
	Let $w$ and $w_n$ be as defined in the preparatory step above.
	Interpret the $\L^1(\Omega_{\lambda})$-functions $w$ and $w_n$ as absolutely continuous measures $w(x) \cdot \Lm^2$ and $w_n(x) \cdot \Lm^2$
	on $\Omega_{\lambda}$. Then, the sequence of measures
	$w_n(x) \cdot \Lm^2$ weakly* converges to $w(x) \cdot \Lm^2$:
	\[
		w_n(x) \cdot \Lm^2 \weaksto w(x) \cdot \Lm^2 \quad, \qquad \mbox{as} \quad n \to \infty \;.
	\]
\end{lem}
\medskip

\begin{proof}	
	Let $\varphi \in \C_0(\Omega_{\lambda})$ be a test function. By changing variables it follows that
	\[
		\int\limits_{\Omega_{\lambda}} (w_n - w)(x) \cdot \varphi(x) \; dx 
		= \int\limits_a^b \int\limits_0^{\lambda} \left( (w_n - w) \cdot \varphi \right) \circ \xi(t,s)  \cdot \det D \xi(t,s) \; dt \; ds \;.
	\]
	We use the representation of $(w_n - w) \circ \xi$ according to equation \refEq{eqn:ReqDWWn} and
	study the convergence of the three summands in equation \refEq{eqn:ReqDWWn} separately.
	The first summand is estimated by
	\begin{align*}
		\left| \int\limits_a^b \int\limits_0^{\lambda} ( (u_{2,n} - u_2) \cdot \varphi) \circ \xi(t,s) \cdot \det D \xi(t,s) \; dt \; ds \right|
		\leq & \; \|u_{2,n} - u_2\|_{L^1(\Omega_{\lambda})} \|\varphi\|_{\infty} \,,
	\end{align*}
	and the right hand side tends to zero, because the sequence $u_{2,n}$ strictly tends to $u_2$ in $\BV(\Omega)$. 
	For the second summand we write 
	\begin{align*}
		 &\left| \int\limits_a^b \int\limits_0^{\lambda}  (\gamma^*u_{2,n}(s) - \gamma^*u_{0,2}(s)) \cdot \varphi \circ \xi(t,s) \cdot \det D \xi(t,s) \; dt \; ds \right| \\
		&=  \left| \int\limits_a^b (\gamma^*u_{2,n}(s) - \gamma^*u_{0,2}(s))  \cdot  |\gamma'(s)| \Biggl( \frac{\int_0^{\lambda} \varphi \circ \xi(t,s) \cdot \det D \xi(t,s) \; dt}{|\gamma'(s)|} \Biggr) \; ds \right| \,.
	\end{align*}
	Let $k_{\lambda}$ and $K_{\lambda}$ be the bounds on the determinant as in lemma \ref{Lem:UniformDetBound}. By the definition of $k_{\lambda}$, we have
	\[
		k_{\lambda} \leq \det D \xi(0,s) \leq \frac{|\gamma'(s)|}{\beta \cdot m_0} \quad \Leftrightarrow  \quad \frac{1}{|\gamma'(s)|} \leq \frac{1}{\beta \cdot m_0 \cdot k_{\lambda}} \;.
	\]
	And consequently,
	\[
		 \frac{\int\limits_0^{\lambda} \varphi \circ \xi(t,s) \cdot \det D \xi(t,s) \; dt}{|\gamma'(s)|} \leq \lambda \cdot \frac{K_{\lambda}}{\beta \cdot m_0 \cdot k_{\lambda}} \cdot \|\varphi\|_{\infty} \;.
	\]
	By the last result, we further estimate:
	\begin{align*}	
		& \leq \; \lambda \cdot \frac{K_{\lambda}}{\beta \cdot m_0 \cdot k_{\lambda}} \cdot \|\varphi\|_{\infty} \cdot \int\limits_a^b |\gamma^*(u_{2,n}-u_{0,2})(s)| \cdot  |\gamma'(s)| \; ds  \\
		& \leq \; \lambda \cdot \frac{K_{\lambda}}{\beta \cdot m_0 \cdot k_{\lambda}} \cdot \|\varphi\|_{\infty} \cdot \int\limits_{\bd \Omega} \left| \; (u_{2,n} - u_2)|_{\bd \Omega} (x) \; \right| \; d\Hm^1(x)   \\
		& = \; \lambda \cdot \frac{K_{\lambda}}{\beta \cdot m_0 \cdot k_{\lambda}} \cdot \|\varphi\|_{\infty} \cdot  \|(u_{2,n} - u_2)|_{\bd \Omega} \|_{L^1(\bd \Omega, \Hm^1)}  \; .
	\end{align*}
	In the last factor we apply the trace operator for $\BV$-functions
	\begin{align*}
		.|_{\bd \Omega} : \BV(\Omega) & \to L^1(\bd \Omega, \Hm^1) \;, & v & \to v|_{\bd \Omega} \;,
	\end{align*}
	which, according to \cite[theorem 3.88]{AmbrosioBV}, is continuous w.r.t. the strict topology on $\BV(\Omega)$. Hence, the factor  $\|(u_{2,n} - u_2)|_{\bd \Omega} \|_{L^1(\bd \Omega, \Hm^1)} $ 
	also tends to zero as $n$ tends to infinity.
	
	Let $\psi(t,s) := \varphi \circ \xi(t,s) \cdot \det D \xi(t,s)$. Then, by changing the order of integration, we get for the third summand
	\begin{align*}
		\int\limits_a^b \int\limits_0^{\lambda} \int\limits_0^t  & \left( \SP< {c_2^0},{\nabla u_{2,n}} > - f_2^0\right) \circ \xi(\tau,s) \; d\tau \cdot \psi(t,s) \; dt \; ds \\
		& = \;  \int\limits_a^b \int\limits_0^{\lambda} \left( \SP< {c_2^0},{\nabla u_{2,n}} > - f_2^0\right) \circ \xi(\tau,s) \left( \int\limits_{\tau}^{\lambda} \psi(t,s) \;dt \right) \; d\tau \; ds \;.
	\end{align*}
	By the definition of $\psi$ and since $\xi$ is a diffeomorphism, there is a continuous function $h \in \C(\cl{\Omega}_{\lambda})$ such that
	\[
		h \circ \xi(\tau,s) = (\det D \xi(\tau,s))^{-1} \cdot \left(\int\limits_{\tau}^{\lambda} \psi(t,s) \;dt \right) \;.
	\]
	With $h$, we write the last integral as
	\begin{align*}
		& = \;  \int\limits_a^b \int\limits_0^{\lambda} \left( \SP< {c_2^0},{\nabla u_{2,n}} > - f_2^0\right) \circ \xi(\tau,s)  \cdot h \circ \xi(\tau,s) \cdot \det D \xi(\tau,s) \; d\tau \; ds \\
		& = \; \int\limits_{\Omega_{\lambda}} \left( \SP< {c_2^0},{\nabla u_{2,n}} > - f_2^0\right)(x) \cdot h(x) \; dx \;.
	\end{align*}
	Next, we use the fact that $u_2$ solves the PDE $\SP< {c_2^0(x)} , {Du_2}> = f_0^2(x) \cdot \Lm^2$
	to formulate the last result as 
	\begin{align*}
		& = \;  \int\limits_{\Omega_{\lambda}}  \SP< {h(x) \cdot c_2^0(x)},{\nabla u_{2,n}(x)} > \; dx -  \int\limits_{\Omega_{\lambda}} \SP< {h(x) \cdot c_2^0(x)}, { dDu_2(x)} >  \\
		& = \;  \int\limits_{\Omega_{\lambda}}  \SP< {\hat{\varphi}(x)},{\nabla u_{2,n}(x)} > \; dx -  \int\limits_{\Omega_{\lambda}} \SP< {\hat{\varphi}(x)}, { dDu_2(x)} >  \;.
	\end{align*}
	In the second equation we have set $\hat{\varphi}(x) := h(x) \cdot c_2^0(x)$
	as a new test function which belongs to $\C(\cl{\Omega_{\lambda}})^2$.
	
	Owing again to the strict convergence of $u_{2,n}$ to $u_2$, we argue by \cite[proposition 3.15]{AmbrosioBV} that the last integral expression tends to zero
	as $n \to \infty$. Summarizing the three steps above we obtain
	\[
		\int\limits_{\Omega_{\lambda}} (w_n - w)(x) \cdot \varphi(x) \; dx \to 0 \qquad \forall \; \varphi \in \C_0(\Omega_{\lambda}) \;,
	\]
	which means $w_n(x) \cdot \Lm^2 \weaksto w(x) \cdot \Lm^2$ on $\Omega_{\lambda}$. \hfill
\end{proof}
\medskip

Based on the properties of the sequence $w_n$ we obtain an estimate of the difference $w = u_1-u_2$ in the next lemma.
Later on, this estimate will account for the operator $U_\lambda$ to be Lipschitz.
\medskip
	
\begin{lem}\label{Lem:PrepULipschitz}
	The difference $w = u_1-u_2$ satisfies
	\begin{align*}
		\|w & \|_{\L^1(\Omega_{\lambda + h})} \leq  \; (\lambda+h) \cdot C_{\lambda+h} \cdot  \|u_{0,1}-u_{0,2}\|_{\L^1(\bd \Omega, \Hm^1)} \\
		& \; + C_{\lambda+h} \cdot \Lm^2(\Omega) \cdot \left( (\lambda+h) \cdot  \|f_1-f_2\|_{\L^{\infty}(\Omega_\lambda)} + h \cdot  \|f_1-f_2\|_{\L^{\infty}(\Omega_{\lambda+h})} \right) \\
		& \; + C_{\lambda+h} \cdot M_{**} \cdot \left( (\lambda+h) \cdot  \|c_1-c_2\|_{\L^{\infty}(\Omega_\lambda)} + h \cdot  \|c_1-c_2\|_{\L^{\infty}(\Omega_{\lambda+h})} \right) \;.
	\end{align*}
	Here, the factor $C_{\lambda} := \frac{K_{\lambda}}{\beta \cdot m_0 \cdot k_{\lambda}}$
	is an increasing function of $\lambda$.
\end{lem}
\medskip

\begin{proof}
We use the approximation of $w$ by $w_n$ again. Because of the weak* convergence according to lemma \ref{Lem:ConvWWn} and because of the lower semi-continuity of the total variation
w.r.t. the weak* convergence (for the semi-continuity of norms, e.g., see \cite{Aliprantis:InfDimAna}), we have
\[
	\|w\|_{L^1(\Omega_{\lambda})} = \left| w\cdot \Lm^2 \right| (\Omega_{\lambda}) 
	\leq \liminf\limits_{n \to \infty } \left| w_n \cdot \Lm^2 \right| (\Omega_{\lambda}) 
	= \liminf\limits_{n \to \infty } \|w_n\|_{L^1(\Omega_{\lambda})} \;,
\]
and thus we can estimate $\|w_n\|_{L^1(\Omega_{\lambda})}$ instead. 
Using the representation of $w_n$ by equation \refEq{eqn:RepWn}, we obtain
\begin{align*}
	\|w_n & \|_{L^1(\Omega_{\lambda+h})} =  \; \int\limits_a^b \int\limits_0^{\lambda+h} |w_n| \circ \xi(t,s) \cdot \det D \xi(t,s) \; dt \; ds \\
	\leq & \; \int\limits_a^b \int\limits_0^{\lambda+h} |\gamma^*(u_{0,1}-u_{0,2})(s)|  \cdot \det D \xi(t,s) \; dt \; ds \\
	&\; + \int\limits_a^b \int\limits_0^{\lambda+h} \int\limits_0^t |f_{1,0}-f_2^0| \circ \xi(\tau,s) \; d\tau  \cdot \det D \xi(t,s) \; dt \; ds \\
	&\; + \int\limits_a^b \int\limits_0^{\lambda+h} \int\limits_0^t \left| \SP< {(c_{1,0}-c_2^0)}, {\nabla u_{2,n}}>\right| \circ \xi(\tau,s) \; d\tau  \cdot \det D \xi(t,s) \; dt \; ds \;.
\end{align*}
By arguing the same way as in the proof of lemma \ref{Lem:ConvWWn} for the first summand we get
\begin{align*}
	\int\limits_a^b \int\limits_0^{\lambda+h} |\gamma^*(u_{0,1}-u_{0,2})(s)|  \cdot & \det D \xi(t,s) \; dt \; ds 
		\leq  (\lambda+h) \cdot C_{\lambda+h} \cdot  \|u_{0,1}-u_{0,2}\|_{\L^1(\bd \Omega, \Hm^1)} \;.
\end{align*}
For the last summand, let $g(\tau,s,t) := (| c_{1,0}-c_2^0|  \cdot  |\nabla u_{2,n}|) \circ \xi(\tau,s) \cdot \det D \xi(t,s)$.
Then, we estimate
\begin{align*}
	\int\limits_a^b \int\limits_0^{\lambda+h} & \int\limits_0^t | \SP< {(c_{1,0}-c_2^0)}, {\nabla u_{2,n}}>| \circ \xi(\tau,s) \; d\tau  \cdot \det D \xi(t,s) \; dt \; ds \\
	\leq & \; \int\limits_a^b \int\limits_0^{\lambda+h} \int\limits_0^t g(\tau,s,t) \; d\tau \; dt \; ds = \int\limits_a^b \int\limits_0^{\lambda+h} \left( \int\limits_{\tau}^{\lambda+h} g(\tau,s,t)  \; dt \right)\; d\tau \; ds \\
	\leq & \; \int\limits_a^b \int\limits_0^{\lambda} \left( \int\limits_{0}^{\lambda+h} g(\tau,s,t)  \; dt \right)\; d\tau \; ds 
	+ \int\limits_a^b \int\limits_{\lambda}^{\lambda+h} \left( \int\limits_{\lambda}^{\lambda+h} g(\tau,s,t)  \; dt \right)\; d\tau \; ds \;.
\end{align*}
For the inner integrals, we have
\begin{align*}
	\int\limits_{\tau}^{\lambda+h} g( & \tau,s,t)  \; dt = ( | c_{1,0}-c_2^0|  \cdot  |\nabla u_{2,n}|) \circ \xi(\tau,s) \int\limits_{\tau}^{\lambda+h} \det D \xi(t,s)  \; dt \\
	& \leq (| c_{1,0}-c_2^0|  \cdot  |\nabla u_{2,n}|) \circ \xi(\tau,s) \cdot (\lambda+h - \tau) \cdot K_{\lambda+h} \\
	& \leq (| c_{1,0}-c_2^0|  \cdot  |\nabla u_{2,n}|) \circ \xi(\tau,s) \cdot \det D \xi(\tau,s) \cdot (\lambda+h - \tau) \cdot \frac{K_{\lambda+h}}{k_{\lambda+h}} \;.
\end{align*}
In the next step we take away the scaling factor, which is in the transport fields and the right-hand sides of the PDE, by $1/\left<c_1,\nabla T_0 \right> \leq 1/(m_0 \cdot \beta)$:
\begin{align*}
	& \leq (| c_{1}-c_2|  \cdot  |\nabla u_{2,n}|) \circ \xi(\tau,s) \cdot \det D \xi(\tau,s) \cdot (\lambda+h - \tau) \cdot \frac{K_{\lambda+h}}{\beta \cdot m_0 \cdot k_{\lambda+h}} \\
	& = (| c_{1}-c_2|  \cdot  |\nabla u_{2,n}|) \circ \xi(\tau,s) \cdot \det D \xi(\tau,s) \cdot (\lambda+h - \tau) \cdot C_{\lambda+h} \;.
\end{align*}
By the last result we infer  on the one hand that
\begin{align*}
	\int\limits_a^b \int\limits_0^{\lambda}   \int\limits_{0}^{\lambda+h}  g(\tau,s,t)   \; dt  \; d\tau \; ds
	& \leq (\lambda+h) \cdot C_{\lambda+h} \int\limits_{\Omega_{\lambda}} | c_{1}-c_2|(x)  \cdot  |\nabla u_{2,n}|(x) \; dx \\
	& \leq (\lambda+h) \cdot C_{\lambda+h} \| c_{1}-c_2\|_{ \L^{\infty}(\Omega_{\lambda}) }  \cdot  \|\nabla u_{2,n}\|_{\L^1(\Omega_{\lambda})} \;,
\end{align*}
and on the other hand that
\begin{align*}
	\int\limits_a^b \int\limits_{\lambda}^{\lambda+h} \int\limits_{\lambda}^{\lambda+h} g(\tau,s,t)  \; dt \; d\tau \; ds 
	& \leq h \cdot C_{\lambda+h} \int\limits_{\Omega_{\lambda+h}} | c_{1}-c_2|(x)  \cdot  |\nabla u_{2,n}|(x) \; dx \\
	& \leq h \cdot C_{\lambda+h} \| c_{1}-c_2\|_{ \L^{\infty}(\Omega_{\lambda+h}) }  \cdot  \|\nabla u_{2,n}\|_{\L^1(\Omega_{\lambda+h})} \;.
\end{align*}
Finally, for the summand 
\[
	\int\limits_a^b \int\limits_0^{\lambda+h} \int\limits_0^t |f_{1,0}-f_2^0| \circ \xi(\tau,s) \; d\tau  \cdot \det D \xi(t,s) \; dt \; ds \;,
\]
we need to perform the same steps with
\[
	g(\tau,s,t) :=( |f_{1,0}-f_2^0|   \cdot  \1_{\Omega_{\lambda+h}} ) \circ \xi(\tau,s) \cdot \det D \xi(t,s) \;,
\]
and end up with the same estimates, but  $\| c_{1}-c_2\|_{ \L^{\infty}(\Omega_{\lambda+h}) }$ has to be replaced by $\| f_{1}-f_2\|_{ \L^{\infty}(\Omega_{\lambda+h}) }$ and 
$\|\nabla u_{2,n}\|_{\L^1(\Omega_{\lambda+h})}$ has to be replaced by $\Lm^2(\Omega_{\lambda+h})$.

Summarizing the last considerations we have an estimate for $\|w_n  \|_{L^1(\Omega_{\lambda+h})}$ by
\begin{align*}
	& \|w_n  \|_{\L^1(\Omega_{\lambda + h})} \leq  \; (\lambda+h) \, C_{\lambda+h} \cdot  \|u_{0,1}-u_{0,2}\|_{\L^1(\bd \Omega, \Hm^1)} \\
	& \; + C_{\lambda+h} \, \Lm^2(\Omega_{\lambda+h}) \left( (\lambda+h) \cdot  \|f_1-f_2\|_{\L^{\infty}(\Omega_\lambda)} + h \cdot  \|f_1-f_2\|_{\L^{\infty}(\Omega_{\lambda+h})} \right) \\
	& \; + C_{\lambda+h} \, \|\nabla u_{2,n}\|_{\L^1(\Omega_{\lambda+h})} \left( (\lambda+h) \cdot  \|c_1-c_2\|_{\L^{\infty}(\Omega_\lambda)} + h \cdot  \|c_1-c_2\|_{\L^{\infty}(\Omega_{\lambda+h})} \right).
\end{align*}
Because $u_{2,n}$ strictly tends to $u_2$, we have
\[
	\|\nabla u_{2,n}\|_{\L^1(\Omega_{\lambda})}  \to |Du_2|(\Omega_{\lambda})  \;.
\]
Hence, going over to the $\liminf$ and plugging in the estimates
\begin{align*}
	|Du_2|(\Omega_{\lambda}) & \leq M_{**}  \qquad\quad \mbox{and} & \Lm^2(\Omega_{\lambda+h}) &\leq  \Lm^2(\Omega) 
\end{align*}
finally yields the assertion.

Because  by lemma  \ref{Lem:UniformDetBound} we know that $K_{\lambda}$ increases while $k_{\lambda}$ decreases with $\lambda$, it
is clear that $C_{\lambda}$ increases with $\lambda$. \hfill
\end{proof}
\medskip

As a corollary, we obtain that $U_\lambda$ is Lipschitz.
\medskip

\begin{cor}\label{Cor:ULipschitz}
	For any choice of $0 < \lambda < 1$ the operator $U_\lambda:\FX_{\lambda} \to \FX_{\lambda}$ is $\L^1$-Lipschitz
	\[
		\|U[v_1]-U[v_2]\|_{\L^1(\Omega_{\lambda})} \leq \lambda \cdot \kappa_{\lambda} \cdot \|v_1-v_2\|_{\L^1(\Omega_{\lambda})} \;.
	\]
	Here, $\kappa_{\lambda}$ is defined by 
	\[
		\kappa_{\lambda} := C_{\lambda} \cdot (L_2 \cdot \Lm^2(\Omega) + L_1 \cdot M_{**}) \;,
	\]
	and is an increasing function of $\lambda$.
\end{cor}
\medskip

\begin{proof}
	Let $v_1, v_2 \in \FX_\lambda$.
	When we consider the operator $U$ or $U_\lambda$, we always use the same boundary data $u_0 \in \FB$. Hence, as mentioned in equation  \refEq{eqn:PrepareU} of the preparatory part at the beginning of this section,
	we have $w = u_1 - u_2 = U_\lambda[v_1] - U_\lambda[v_2]$, since $u_{0,1}=u_{0,2}=u_0$. By using lemma \ref{Lem:PrepULipschitz} with $h=0$ we see that
	\begin{align*}
		\|U[v_1]-U[v_2]\|_{\L^1(\Omega_{\lambda})} \leq
		  \lambda \cdot C_{\lambda} \cdot \Bigl(  \Lm^2(\Omega) &\cdot  \|f_1-f_2\|_{\L^{\infty}(\Omega_\lambda)}  
		+ M_{**} \cdot  \|c_1-c_2\|_{\L^{\infty}(\Omega_\lambda)}  \Bigr).
	\end{align*}
	For the differences $f_1-f_2$ and $c_1-c_2$ we use now the functional causality and the Lipschitz conditions, which we require. That is
	\begin{align*}
		\|f_1-f_2\|_{\L^{\infty}(\Omega_{\lambda})} &= \|f[v_1]-f[v_2]\|_{\L^{\infty}(\Omega_\lambda)} = \| (f[v_1]-f[v_2]) |_{\Omega_\lambda}\|_{\infty}  \\
		& = \| (f[v_1\cdot \1_{\Omega_{\lambda}}]-f[v_2 \cdot \1_{\Omega_\lambda}]) |_{\Omega_\lambda}\|_{\infty} \leq \| f[v_1\cdot \1_{\Omega_{\lambda}}]-f[v_2 \cdot \1_{\Omega_{\lambda}}] \|_{\infty} \\
		& \leq L_2 \cdot \| v_1\cdot \1_{\Omega_{\lambda}}-v_2 \cdot \1_{\Omega_{\lambda}}\|_{\L^1(\Omega)} = L_2 \cdot \| v_1-v_2 \|_{\L^1(\Omega_{\lambda})} 
	\end{align*}
	and analogously
	\[
		\|c_1-c_2\|_{\L^{\infty}(\Omega_\lambda)} \leq 	L_1 \cdot \| v_1-v_2 \|_{\L^1(\Omega_{\lambda})} \;.
	\]
	Putting everything together, it follows that
	\begin{align*}
		\|U[v_1]-U[v_2]\|_{\L^1(\Omega_{\lambda})} &\leq \lambda \cdot C_{\lambda} \cdot (L_2 \cdot \Lm^2(\Omega) + L_1 \cdot M_{**}) \cdot \| v_1-v_2 \|_{\L^1(\Omega_{\lambda})} \\
		& \leq \lambda \cdot \kappa_{\lambda} \cdot \|v_1-v_2\|_{\L^1(\Omega_{\lambda})} \;.
	\end{align*}
	Finally, $\kappa_{\lambda}$ increases with $\lambda$, since $C_{\lambda}$ does so, too. \hfill
\end{proof}
\medskip

Now that we have brought together all ingredients we are able to show the uniqueness of the fixed point.
\medskip

\begin{theo}\label{Theo:UniquenessQL}
	(Uniqueness)
	The solution operator $U:\FX \to \FX$ of the (non-restricted) original problem \refEq{eqn:LinOfQuasiLin} has a unique fixed point $u \in \FX$, $u=U[u]$.
\end{theo}
\medskip

\begin{proof}
	First, we show that, for any choice of $0 < \lambda < 1$, the  domain-restricted operator $U_\lambda:\FX_{\lambda} \to \FX_{\lambda}$ has a unique fixed point. 
	In order to do so we decompose $\Omega_{\lambda}$ into finitely many stripes $\Omega_{(l+1)h,lh}$ of ''thickness'' $h >0 $.
	For such stripes, we use the notation
	\[
		\Omega_{\lambda+h,\lambda} = \Omega_{\lambda+h} \wo \Omega_{\lambda} = \{z \in \Omega: \lambda \leq T_0(z) < \lambda + h\} \;.
	\]
	Let the step size $h$ be such that $h < 1/\kappa_{\lambda}$, and let $L = \left\lfloor \lambda/h \right\rfloor \in \N$ be the number of steps. Then,
	\[
		\Omega_{\lambda} = \bigcup\limits_{l=0}^{L-1} \Omega_{(l+1)h,lh}  \cup \Omega_{\lambda,Lh} \;.
	\]
	For the first step, consider the operator $U_h:\FX_{h} \to \FX_{h}$ on $\Omega_h = \Omega_{h,0}$. By corollary \ref{Cor:ULipschitz} and the choice of $h$ we have
	a contraction
	\[
		\|U_h[v_1]-U_h[v_2]\|_{\L^1(\Omega_{h})} \leq h \cdot \kappa_{h} \cdot \|v_1-v_2\|_{\L^1(\Omega_{h})} \leq h \cdot \kappa_{\lambda} \cdot \|v_1-v_2\|_{\L^1(\Omega_{h})} \;.
	\]
	If now $u_1 = U[u_1]$ and $u_2 = U[u_2]$ are two fixed points, we have, after domain-restriction onto  $\Omega_{h}$,
	\[
		\|u_1-u_2\|_{\L^1(\Omega_{h})} \leq h \cdot \kappa_{\lambda} \cdot \|u_1-u_2\|_{\L^1(\Omega_{h})}
	\]
	and consequently
	\[
		0 \leq (1-h \cdot \kappa_{\lambda}) \|u_1-u_2\|_{\L^1(\Omega_{h})} \leq 0 \;.
	\]
	Hence, all fixed points coincide on the stripe $\Omega_{h}$.
		
	Next, we perform an inductive step. Assume that all fixed points coincide on $\Omega_{lh}$, we show that they must also coincide on $\Omega_{(l+1)h} = \Omega_{lh + h}$.
	Let $u_1 = U[u_1]$ and $u_2 = U[u_2]$ be two fixed points again. With $w=u_1-u_2$, by lemma \ref{Lem:PrepULipschitz} we know that
	\begin{align*}
		& \|u_1-u_2  \|_{\L^1(\Omega_{(l+1)h})} \leq C_{(l+1)h} \; \cdot\\
		& \; \biggl( \Lm^2(\Omega) \cdot \left( (l+1)h \cdot  \|f[u_1]-f[u_2]\|_{\L^{\infty}(\Omega_{lh})} + h \cdot  \|f[u_1]-f[u_2]\|_{\L^{\infty}(\Omega_{(l+1)h})} \right) \\
		& \; +  M_{**} \cdot \left( (l+1)h \cdot  \|c[u_1]-c[u_2]\|_{\L^{\infty}(\Omega_{lh})} + h \cdot  \|c[u_1]-c[u_2]\|_{\L^{\infty}(\Omega_{(l+1)h})} \right) \biggr) \;.
	\end{align*}
	Because $u_1$ and $u_2$ coincide on $\Omega_{lh}$, we have
	\[
		\|f[u_1]-f[u_2]\|_{\L^{\infty}(\Omega_{lh})} = 0 \quad \mbox{  and  } \quad \|c[u_1]-c[u_2]\|_{\L^{\infty}(\Omega_{lh})} = 0 \;,
	\]
	and thus, the estimate reduces to
	\begin{align*}
		 \|u_1-u_2  \|_{\L^1(\Omega_{(l+1)h})} \leq 
		& \;  C_{(l+1)h} \cdot h \cdot ( \Lm^2(\Omega)  \cdot  \|f[u_1]-f[u_2]\|_{\L^{\infty}(\Omega_{(l+1)h})} \\
		& \qquad + M_{**}  \cdot  \|c[u_1]-c[u_2]\|_{\L^{\infty}(\Omega_{(l+1)h})} )\;.
	\end{align*}
	By using the Lipschitz conditions on $c$ and $f$ again we have
	\begin{align*}
		\|u_1-u_2  \|_{\L^1(\Omega_{(l+1)h})} &\leq h \kappa_{(l+1)k} \cdot \|u_1-u_2  \|_{\L^1(\Omega_{(l+1)h})} \\
                &\leq h \kappa_{\lambda} \cdot \|u_1-u_2  \|_{\L^1(\Omega_{(l+1)h})} \,.
	\end{align*}
	More precisely, because we have assumed that $u_1$ and $u_2$ coincide on $\Omega_{lh}$, the latter inequality in fact means
	\[
		\|u_1-u_2  \|_{\L^1(\Omega_{(l+1)h,lh})}  \leq h \kappa_{\lambda} \cdot \|u_1-u_2  \|_{\L^1(\Omega_{(l+1)h,lh})} \;.
	\]
	By the contractiveness, $h \kappa_{\lambda} < 1$, we see  that $\|u_1-u_2  \|_{\L^1(\Omega_{(l+1)h,lh})} =0$ and so the fixed points also coincide on the next stripe $\Omega_{(l+1)h,lh}$.
		
	For the last stripe we have to adapt the step size to $\hat{h} = \lambda - Lh$.
	But, since $\hat{h} \leq h$, the same argumentation applies.
		
	As claimed before, the  domain restricted operator $U_\lambda:\FX_{\lambda} \to \FX_{\lambda}$ has a unique fixed point for any choice of $0 < \lambda < 1$.
	Now the last step: assume by contradiction that the non-restricted operator $U:\FX \to \FX$ has two different fixed points, $u_1$ and $u_2$. Therefor, $u_1$ and $u_2$ must
	differ on a subset $W \subset \Omega$ with $\Lm^2(W) \neq 0$. Because the stop set $\Sigma$ has Lebesgue measure zero, $\Lm^2(\Sigma)= 0$, we can choose
	$0 < \lambda < 1$ so close to 1 that $\Lm^2(\Omega_{\lambda} \cap W) \neq 0$. Thus, we have
	\[
		\|u_1-u_2  \|_{\L^1(\Omega_{\lambda})} = \|u_1-u_2  \|_{\L^1(\Omega_{\lambda}\cap W)} \neq 0 \;.
	\]
	But, because $u_1|_{\Omega_{\lambda}}$ and $u_2|_{\Omega_{\lambda}}$ are
	fixed points of the  domain restricted operator $U_\lambda:\FX_{\lambda} \to \FX_{\lambda}$, we also have
	\[
		\|u_1-u_2  \|_{\L^1(\Omega_{\lambda})} = 0
	\]
	by the previous uniqueness proof. A contradiction. \hfill
\end{proof}

\section{Continuous Dependence of the Fixed Point}\label{Sect:ContDep}
In this section, we will show that the unique fixed point depends $\L^1$-con\-tinuously on the following data: the transport field, the right-hand side, and the boundary data.
We consider two linear problems:
\begin{align*}
	\SP< {c[v](x)} , {Du} > &= f[v](x) \cdot \Lm^2 \quad \mbox{ in } \quad \Omega \wo \Sigma \;, \\
	u|_{\bd \Omega} &= u_0 \;,
\end{align*}  
and 
\begin{align*}
	\SP< {\tilde{c}[\tilde{v}](x)} , {D\tilde{u}} > &= \tilde{f}[\tilde{v}](x) \cdot \Lm^2 \quad \mbox{ in } \quad \Omega \wo \Sigma \;, \\
	\tilde{u}|_{\bd \Omega} &= \tilde{u}_0 \;,
\end{align*} 
where we assume that for both problems the same domain $\Omega$, the same stop set $\Sigma$, and the same time function $T$ (with transformed version $T_0$) are specified.

Moreover, we assume that $c$ and $\tilde{c}$ both satisfy the requirement \ref{Req:TransportFieldFunc} with the same bounds,
and that $f$ and $\tilde{f}$ both satisfy the requirement \ref{Req:RHSFunc} with the same bounds. Finally, we assume $u_0 \in \FB$ and $\tilde{u}_0 \in \FB$.
By the latter assumption we are sure that we obtain two solution operators 
\begin{align*}
	U:\FX & \to \FX \;, & v & \to U[v] \;, \\
	\tilde{U}:\FX & \to \FX \;, & \tilde{v} & \to \tilde{U}[\tilde{v}] \;, 
\end{align*}
which respectively correspond to the two linear problems above and possess the same domain and range $\FX$, which depends on all those bounds.

We view $\tilde{c}$ and $\tilde{f}$ as perturbed versions of $c$ and $f$. In order to measure the perturbation we introduce the following norm:
\medskip

\begin{defi}
	For  maps of type $g:\L^1(\Omega) \to \C_b(\Omega \wo \Sigma)^d$ or of type $g:\L^1(\Omega) \to \C(\cl{\Omega})^d$, $d \in \N$, we define the norm
	\[
		\|g\|_0 := \sup\limits_{v \in \L^1(\Omega)} \|g[v]\|_{\infty} \;.
	\]
\end{defi}
\medskip

\begin{theo}\label{Theo:ContDepQL}
	(Continuous dependence)
	Consider the two solution operators $U$ and $\tilde{U}$ as described above. Let $u=U[u]$ and $\tilde{u}=\tilde{U}[ \tilde{u} ]$ be the unique fixed points of these operators.
	\smallskip
	
	Then, for every $\varepsilon > 0 $, one can find $\delta > 0$ such that $\|u-\tilde{u}\|_{\L^1(\Omega)} \leq \varepsilon$,
	whenever
	\[
		\left( \|u_0-\tilde{u}_0\|_{\L^1(\bd \Omega, \Hm^1)} \; + \; \Lm^2(\Omega) \cdot \|f-\tilde{f}\|_0 \; + \; M_{**} \cdot \|c-\tilde{c}\|_0 \right) \; \leq \; \delta \;.
	\]
\end{theo}
\medskip

\begin{proof}
	Let $v, \tilde{v} \in \FX$ be arbitrary but fixed, and set $u_1 := U[v]$ and $u_2 := \tilde{U}[\tilde{v}]$.
	As before, we set up a PDE for the difference $w := u_1 -u_2$ on the restricted domain $\Omega_{\lambda}$, $0 < \lambda < 1$ by
	\begin{align*}
		\SP< {c[v](x)}, {Dw} > &= (f[v]-\tilde{f}[\tilde{v}])(x) \cdot \Lm^2 - \SP< {(c[v]-\tilde{c}[\tilde{v}])(x)}, {Du_2}> \quad \mbox{ in } \quad \Omega_{\lambda} \;, \\
		w|_{\bd \Omega} &= u_0 - \tilde{u}_0 \;.
	\end{align*}
	Again, we choose a sequence $u_{2,n} \in \C^{\infty}(\Omega)$ which strictly converges to $u_2$ in $\BV(\Omega)$.
	And again, the sequence $w_n$ of solutions to the approximate PDE -- which has the absolutely continuous measure $\nabla u_{2,n}(x) \cdot \Lm^2$ instead of $Du_2$ -- converges weakly* to $w$.
		
	In order to proceed as in lemma \ref{Lem:PrepULipschitz} we supplement the right-hand side of the approximate PDE to
	\begin{align*}
		\SP< {c[v](x)}, {Dw_n} > =& \; \left( (f[v]-f[\tilde{v}])  - \SP< {c[v]-c[\tilde{v}]}, {\nabla u_{2,n}}> \right)(x) \cdot \Lm^2 \\
		& \; + \left( (f[\tilde{v}]-\tilde{f}[\tilde{v}]) - \SP< {c[\tilde{v}]-\tilde{c}[\tilde{v}]}, {\nabla u_{2,n}}> \right)(x) \cdot \Lm^2 \;.
	\end{align*}
	For the first summand of the new right-hand side we will apply the steps from the proof of lemma \ref{Lem:PrepULipschitz}.
	For the second summand we operate in a simpler way. Let $\xi = \xi[v]$ be the characteristics corresponding to the field $c_0[v]$ and let
	\[
		g = \frac{\left| (f[\tilde{v}]-\tilde{f}[\tilde{v}]) - \SP< {c[\tilde{v}]-\tilde{c}[\tilde{v}]}, {\nabla u_{2,n}}> \right| }{\SP<{c[v]},{\nabla T_0}>} \;.
	\]
	As in the proof of lemma \ref{Lem:PrepULipschitz} we have to estimate
	\[
		\int\limits_a^b \int\limits_0^{\lambda+h} \int\limits_0^t g \circ \xi(\tau,s) \; d\tau \det D\xi(t,s) \; dt \; ds \leq \ldots \;.
	\]
	After having changed the order of integration and having estimated the determinant, we obtain
	\begin{align*}
		\ldots & \leq (\lambda + h) \cdot \frac{K_{\lambda+h}}{k_{\lambda+h}} \cdot  \int\limits_a^b \int\limits_0^{\lambda+h} g \circ \xi(\tau,s)  \det D\xi(\tau,s) \; d\tau \; ds \\
		\leq & (\lambda + h) \cdot  C_{\lambda+h} \cdot \int\limits_{\Omega_{\lambda+h}} \left| (f[\tilde{v}]-\tilde{f}[\tilde{v}]) - \SP< {c[\tilde{v}]-\tilde{c}[\tilde{v}]}, {\nabla u_{2,n}}> \right| (x) \; dx \\
		\leq & (\lambda + h) \cdot  C_{\lambda+h} \left( \Lm^2(\Omega) \cdot \| f[\tilde{v}]-\tilde{f}[\tilde{v}] \|_{\infty} + \|\nabla u_{2,n}\|_{\L^1(\Omega)} \cdot \| c[\tilde{v}]-\tilde{c}[\tilde{v}] \|_{\infty} \right) \\
		\leq & (\lambda + h) \cdot  C_{\lambda+h} \left( \Lm^2(\Omega) \cdot \| f-\tilde{f} \|_{0} + \|\nabla u_{2,n}\|_{\L^1(\Omega)} \cdot \| c-\tilde{c} \|_{0} \right) \;.
	\end{align*}
	Putting both estimates together and then going over to the $\liminf$, we end up with
	\begin{align*}
		& \|w\|_{\L^1(\Omega_{\lambda + h})} \leq  \; (\lambda+h) \cdot C_{\lambda+h} \cdot  \|u_{0}-\tilde{u}_{0}\|_{\L^1(\bd \Omega, \Hm^1)} \\
		& \; + C_{\lambda+h} \cdot \Lm^2(\Omega) \left( (\lambda+h) \cdot  \|f[v]-f[\tilde{v}]\|_{\L^{\infty}(\Omega_\lambda)} + h \cdot  \|f[v]-f[\tilde{v}]\|_{\L^{\infty}(\Omega_{\lambda+h})} \right) \\
		& \; + C_{\lambda+h} \cdot M_{**} \left( (\lambda+h) \cdot  \|c[v]-c[\tilde{v}]\|_{\L^{\infty}(\Omega_\lambda)} + h \cdot  \|c[v]-c[\tilde{v}]\|_{\L^{\infty}(\Omega_{\lambda+h})} \right) \\
		& \; + (\lambda + h) \cdot  C_{\lambda+h} \left( \Lm^2(\Omega) \cdot \| f-\tilde{f} \|_{0} + M_{**} \cdot \| c-\tilde{c} \|_{0} \right) \;.
	\end{align*}
	
	Now we can show the continuous dependence in the domain restricted situation. Fix $0 < \lambda < 1$, choose a step size $0 < h < 1/\kappa_{\lambda} $ and let
	$L = \left\lfloor \lambda/h \right\rfloor \in \N$ be the number of steps. Furthermore, let 
	\[
		\left( \|u_0-\tilde{u}_0\|_{\L^1(\bd \Omega, \Hm^1)} \; + \; \Lm^2(\Omega) \cdot \|f-\tilde{f}\|_0 \; + \; M_{**} \cdot \|c-\tilde{c}\|_0 \right) \; \leq \; \delta \;,
	\]
	for some $\delta > 0$. Let $l \in \N_0$, $l \leq L$.  With the result above we estimate on the set $\Omega_{(l+1)h}$:
	\begin{align*}
		\|w & \|_{\L^1(\Omega_{(l+1) h})} \leq  \; \lambda\cdot C_{\lambda} \cdot \delta \\
		& \; + C_{\lambda} \cdot \Lm^2(\Omega) \cdot \left( \lambda \cdot  \|f[v]-f[\tilde{v}]\|_{\L^{\infty}(\Omega_{lh})} + h \cdot  \|f[v]-f[\tilde{v}]\|_{\L^{\infty}(\Omega_{(l+1)h})} \right) \\
		& \; + C_{\lambda} \cdot M_{**} \cdot \left( \lambda \cdot  \|c[v]-c[\tilde{v}]\|_{\L^{\infty}(\Omega_{lh})} + h \cdot  \|c[v]-c[\tilde{v}]\|_{\L^{\infty}(\Omega_{(l+1)h})} \right) \;.
	\end{align*}
	By using the functional causality, the Lipschitz conditions on $c$ and $f$, and the definition of $\kappa_{\lambda}$ from corollary \ref{Cor:ULipschitz} we obtain
	\begin{align*}
		\|w & \|_{\L^1(\Omega_{(l+1) h})} \leq  \; \lambda \cdot C_{\lambda} \cdot \delta + \lambda \kappa_{\lambda} \cdot \|v-\tilde{v}\|_{\L^{1}(\Omega_{lh})}  
		+ h \kappa_{\lambda} \cdot \|v-\tilde{v}\|_{\L^{1}(\Omega_{(l+1)h})}  \,.
	\end{align*}
	
	Let $\hat{\delta} = \lambda \cdot C_{\lambda} \cdot \delta$. Now, we plug in the two fixed points $u$ and $\tilde{u}$, i.e., we set $u_1 = v = u$ and $u_2 = \tilde{v} = \tilde{u}$,
	\begin{align*}
		\|u-\tilde{u}\|_{\L^{1}(\Omega_{(l+1)h})} \leq & \;  \hat{\delta} + \lambda \cdot \kappa_{\lambda} \cdot \|u-\tilde{u}\|_{\L^{1}(\Omega_{lh})}  
		+ h \cdot \kappa_{\lambda} \cdot \|u-\tilde{u}\|_{\L^{1}(\Omega_{(l+1)h})}  \;.
 	\end{align*}
	We define the error $e_l$ on the set $\Omega_{lh}$ by $e_l := \|u-\tilde{u}\|_{\L^{1}(\Omega_{lh})}$.
	Then, by our choice of $h$, the last estimate yields the error recursion
	\[
		0 \leq (1-h\kappa_{\lambda}) \cdot e_{l+1} \leq \hat{\delta} + \lambda \cdot \kappa_{\lambda} \cdot e_l \;,
	\]
	which leads to
	\[
		e_{l+1} \leq \sum_{k=0}^l \alpha^k  \cdot \frac{\hat{\delta}}{1-h\kappa_{\lambda}} \quad \mbox{with} \quad \alpha := \frac{\lambda \cdot \kappa_{\lambda}}{1-h\kappa_{\lambda}} \;.
	\]
	In summary, we get
	\[
		\|u-\tilde{u}\|_{\L^{1}(\Omega_{\lambda})} \leq e_{L+1} \leq \left( \frac{1-\alpha^{L+1}}{1-\alpha } \cdot \frac{\lambda \cdot C_{\lambda} }{1-h\kappa_{\lambda}} \right) \cdot \delta 
	\]
	and the continuous dependence for the domain restricted case is obvious.
		
	Final step: let $\varepsilon > 0 $. For the full domain $\Omega$ we choose $\lambda$ so close to 1 that
	\[
		\|u-\tilde{u}\|_{\L^{1}(\Omega \wo \Omega_{\lambda})} \leq \frac{\varepsilon}{2} \;.
	\]
	In dependence of this $\lambda$ we find $h$ and $L$. What remains to do is to require
	\[
		\delta = \left( \frac{1-\alpha^{L+1}}{1-\alpha } \cdot \frac{\lambda \cdot C_{\lambda} }{1-h\kappa_{\lambda}} \right)^{-1} \cdot \frac{\varepsilon}{2} \;,
	\]
	then, we get $\|u-\tilde{u}\|_{\L^1(\Omega)} \leq \varepsilon$, whenever
	\[
		\left( \|u_0-\tilde{u}_0\|_{\L^1(\bd \Omega, \Hm^1)} \; + \; \Lm^2(\Omega) \cdot \|f-\tilde{f}\|_0 \; + \; M_{**} \cdot \|c-\tilde{c}\|_0 \right) \; \leq \; \delta \;.
	\]
	\hfill
\end{proof}

\section{Discussion}
We have shown in this paper that the quasi-linear problem \refEq{eqn:QuasiLinProblem} has a unique solution which depends continuously on the coefficients of the PDE
and the boundary data. The special ingredients here have been the functional causality and the Lipschitz continuity of the coefficients.

The good behavior of the time function $T$ and the transport field $c$ close to and on the stop set $\Sigma$ have not been used in sections \ref{Sect:Unique} and \ref{Sect:ContDep}.
So, in the theory presented we could consider solutions which belong only to $\BV(\Omega \wo \Sigma)$.

In \cite{MyPaper1}, the mentioned good behavior of $T$ and $c$ close to and on $\Sigma$ accounted for closing the gap $\Sigma$: for the linear case we obtained that the solution
belongs to $\BV(\Omega)$ and not only $\BV(\Omega \wo \Sigma)$, see  \cite[theorem 4.1]{MyPaper1}. This feature was crucial to proof the existence of a solution
to the quasi-linear case (even for non-causal functional dependence) by the Schauder fixed point theorem because $\FX \subset \BV(\Omega)$ is now compact.
Moreover, the fixed point solutions belong themselves to $\BV(\Omega)$. So here, we have again the closing the gap feature and additionally uniqueness together with
continuous dependence.   

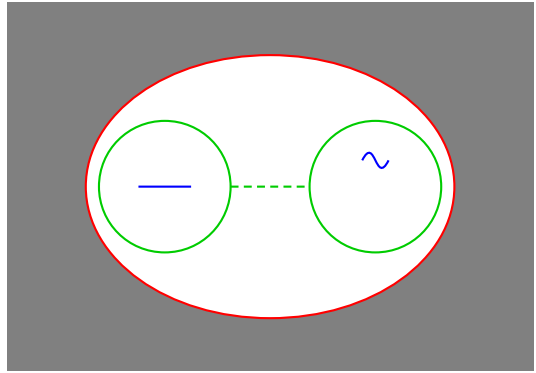
\begin{figure}[!ht]
	\begin{center}
		\begin{tikzpicture}[scale=7]
			\colorlet{darkgreen}{green!80!black}
			
			\fill[gray] (0,0) rectangle (1,0.7);
			\filldraw[fill=white,draw=red,thick] (0.5,0.35) ellipse (0.35cm and  0.25cm);
			
			\draw[darkgreen,thick] (0.3,0.35) circle (0.125cm);
			\draw[darkgreen,thick] (0.7,0.35) circle (0.125cm);
			\draw[darkgreen,thick,densely dashed] (0.425,0.35) -- (0.575,0.35); 
						
			\draw[blue,thick] (0.25,0.35) -- (0.35,0.35);
			\draw[blue,thick] (0.675,0.4) .. controls (0.7,0.45) and (0.7,0.35) .. (0.725,0.4);
		\end{tikzpicture}
	\end{center}
	\caption{Three levels of a time function $T$ with a locally maximal level in the middle (dashed green line). 
		The white area is the domain $\Omega$ with its boundary, which is the start level $T=0$, in red.
		The solid green and the dashed green lines together are a saddle level (the terminal nodes of the dashed green line behave like saddle nodes).
		The dashed green line alone is an intermediate stop set while the blue lines are the maximal level of $T$, i.e., the final stop set which is disjoint.}  
	\label{Fig:Saddle}
\end{figure}

The advantage of closing the gap here is, that we can allow for more general time functions. For example we can use one with levels as illustrated in figure \ref{Fig:Saddle}.
Here in a first step one solves from the start level to the saddle level. With our closing the gap feature we can even reach the saddle level and close the first gap which is given by the intermediate stop set.
Finally, one solves the remaining two problems and gets a unique global solution in $\Omega$. Without closing the gap, one could only produce a local solution which is defined on
the region between the start and the saddle level.  For more on this topic see \cite[chapter 5]{MyDiss}.

\section{Acknowledgements}
The author would like to thank Folkmar Bornemann for his advice and the inspiring discussions.

\bibliographystyle{siam}
\bibliography{DPaper2}

\end{document}